%
%
%
%
\def\res{\restriction}

\let\he=\heads
\let\hee=\headss

\def\hh #1\par{\underbar{\he#1}\vskip1pc}
\def\hhh #1\par{\underbar{\hee#1}\vskip1.2pc}
\def\ar #1\ {\bf#1}
\font\title=cmbx12
\font\by=cmr8
\font\author=cmr10
\font\adress=cmsl10
\font\abstract=cmr8

\def\tit #1\par{\centerline{\title #1}}
\def\bby{\centerline{\by BY}\par}
\def\aut #1\par{\centerline{\author #1}}
\def\abs{\abstract\centerline{ABSTRACT}\par}
\def\abss #1\par{\abstract\midinsert\narrower\narrower\noindent #1\endinsert}
\def\\{\noindent}
\def\vv{\par\rm}
\def\vvv{\medbreak\rm}

\def\sec{\vskip 2pc}

\def\sh{[Sh]}

\def\gu{[GuSh]}

\def\qed{\line{\hfill$\heartsuit$}\medbreak}
\def\pt #1. {\medbreak\noindent{\bf#1. }\enspace\sl}
\def\pd #1. {\medbreak \noindent{\bf#1. }\enspace}
\def\dd{Definition\ } 
\def\bb #1{{\bar #1}}
\def\bbp{\bar P}
\def\fo{\ ||\kern -4pt-}
\def\om{\omega}
\def\al{\alpha}
\def\be{\beta}
\def\ga{\gamma}
\def\Ga{\Gamma}

\def\de{\delta}
\def\et{\eta}
\def\si{\sigma}
\def\ta{\tau}
\def\ale{\aleph_0}
\def\vf{\varphi}
\def\sb{\subseteq}

\def\imp{\Rightarrow}
\def\lan{\langle}
\def\ran{\rangle}
\def\all{\forall}
\def\ex{\exists}
\def\mo{\models}
\def\qmu{$Q$\llap{\lower5pt\hbox{$\sim$}\kern1pt}$_\mu$}

\def\M{{\cal M}}
\def\T{{\cal T}}
\def\em{\emptyset}

\def\sm{\setminus}
\def\and{\ \&\ }
\def\th{{\rm Th}}


\def\st{such that\ }

\def\int{interpretation\ }

\def\today{\ifcase\month\or January\or February\or March\or April\or May
\or June\or July\or August\or September\or October\or November\or December\fi 
\space\number\day, \number\year}
\def\ee{\thinspace^\wedge}

\def\seg{\triangleleft}
\def\bina{\thinspace^{\omega>}2}
%
%
%

\baselineskip 1.2pc
\hsize=6in
\input mssymb  

\def\seg{\triangleleft}
\def\hdc{{\rm Hdeg}$(C)$}
\def\Q{{\cal Q}}
\def\vf{\varphi}                                                       
\def\ep{\epsilon}                                                      
\def\B{{\cal B}}                                                
\def\P{{\cal P}}                                                
\def\bpl{\bigoplus}                                             
\def\ot{{\rm otp}}                                              
\def\cf{{\rm cf}}                                               
\def\lg{{\rm lg}}                                               
\def\dpp{{\rm dp}}                                               
\line{\hfill\today}
\vskip.5in
\tit Uniformization and Skolem Functions in the Class of Trees.
\par
\vskip.5in
\bby
\aut SHMUEL LIFSCHES and SAHARON SHELAH\footnote*{The second author would like 
to thank the U.S.--Israel Binational Science Foundation for 
partially supporting this research. Publ. ***}
\par
{\adress Institute of Mathematics, The Hebrew University of Jerusalem, 
Jerusalem, Israel} \par
\vskip.5in
%
%
%
\abs
\abss The monadic second-order theory of trees allows quantification over 
elements and over arbitrary subsets. We classify the class of trees with 
respect to the question: does a tree $T$ have definable Skolem functions 
(by a monadic formula with parameters)?  This continues [LiSh539] where 
the question was asked only with respect to choice functions. Here 
we define a subclass of the class of tame trees (trees with a definable  
choice function) and prove that this is exactly the class (actually set)
of trees with definable Skolem functions.
\par
\vskip.6in  \rm
%
%
%
\hh 1. Introduction: The Uniformization Problem  \par
%
%
%
\pd \dd 1. The monadic second-order logic is the fragment of the full 
second-order logic that allows quantification over elements and over monadic 
(unary) predicates only. The monadic version of a first-order language $L$ can 
be described as the augmentation of\/ $L$ by a list of quantifiable set 
variables and by new atomic formulas $t\in X$ where $t$ is a first order term 
and $X$ is a set variable. The monadic theory of a structure $\M$ is the 
theory of $\M$ in the extended language where the set variables range over all 
subsets of  $|\M|$ and $\in$ is the membership relation.
\pd \dd 2. The {\sl monadic language of order} $L$ is the monadic version 
of the language of order $\{<\}$. For simplicity, we add to $L$  
the predicate sing$(X)$ saying ``$X$ is a singleton'' and use only formulas 
with set variables.  
Thus the meaning of $X<Y$ is: $X=\{x\}\ \&\ Y=\{y\}\ \&\ x<y$.  \vvv
\pd \dd 3.  Let $T$ be a tree and $\bbp\sb T$.

\\(1) $\vf$ is an $(n,l)$-formula if $\vf=\vf(X,Y,\bbp)$ with  
$\dpp(\vf)=n$ and $l(\bbp)=l$. 

\\(2)$\vf=\vf(X,Y,\bbp)$ is {\sl potentially uniformizable in $T$} (p.u) if 
$T\mo(\all Y)(\ex X)\vf(X,Y,\bbp)$. \vv
\sec
\vfill \eject
%
%
%
\hh 2. Tame Trees  \par
\pd \dd 2.1. A tree is a partially ordered set $(T,\seg)$ \st for every 
$\et\in T$, $\{\nu:\nu\seg\et\}$ is linearily orderd by $\seg$.

\\Note, a chain $(C,<^*)$ and even a set without structure $I$ is a tree. 

\\Branch, Sub-branch, Initial segment. \vv
\pd \dd 2.2.  (1) $(C,<^*)$ is a scattered chain iff ...

\\(2) For a scattered chain $(C,<^*)$ \ \hdc \ is defined inductively by:

\\ \hdc=$0$  iff ...

\\ \hdc=$\al$  iff ...

\\ \hdc$\ge\de$  iff ... \vv
\pt Theorem 2.3. Hdeg$(C)$ exists for every scattered chain $C$. \vv
\pt Lemma 2.4. Hdeg$(C)<\om$ then $C$ has a definable well ordering. 
\vv
\pd Proof. See A1 in the appendix 
\vv \qed
\pd \dd 2.5. $\sim_A^0$, $\sim_A^1$. \ \ \ \ \ \ \ (from  [LiSh539] 4.1)
\vv
\pd \dd 2.6. (1) A tree $T$ is called {\sl wild} if either \par 
$(i)$ $sup\big\{|top(A)/\sim^1_A| : A\sb T\ {\rm an\ initial\ segment}\big\}
\ge\ale$ \ \ or 

$(ii)$ There is a branch $B\sb T$ and an embedding $f\colon\Q\to B$ \ \ or

$(iii)$ All the branches of\/ $T$ are scattered linear orders but  
$sup\big\{{\rm Hdeg}(B) : B\ {\rm a\ branch\ of\ } T\big\} \ge\om$.

$(iv)$ There is an embedding $f\colon\bina\to T$

\\(2) A tree $T$ is {\sl tame for} $(n^*,k^*)$ if 
the value in $(i)$ is $\le n^*$,  $(ii)$ does not hold and the value in 
$(iii)$ is $\le k^*$ 

\\(3) A tree $T$ is {\sl tame} if\/ $T$ is tame for $(n^*,k^*)$ for some 
$n^*,k^*<\om$. \vv
\\

\\The following is the content of [LiSh539], $(2)\imp (3)$ is given in theorem 
A2 in the appendix.
\pt Theorem 2.7. The following are equivalent:

1. $T$ has a definable choice function.

2. $T$ has a definable well ordering.

3. $T$ is tame. 
\vv \qed
\sec
%
%
%
\hh 3. Composition Theorems  \par
%
%
%
\pd Notations.  $x,y,z$ denote individual variables, $X,Y,Z$ are set 
variables, $a,b,c$ elements and $A,B,C$ sets. $\bb a,\ \bb A$  are finite 
sequences and $\lg(\bb a),\ \lg(\bb A)$ their length. We write e.g. 
$\bb a\in C$ and $\bb A\sb C$ instead of $\bb a\in\,^{\lg(\bb a)}C$  or 
$\bb A\in\,^{\lg(\bb A)}{\cal P}(C)$   \vv
\pd \dd 3.1. For any chain $C$, $\bb A\in\,^{\lg(\bb A)}{\cal P}(C)$, 
and a natural number $n$, define by induction 
$$t = \th^n(C;\bb A)$$ 
for $n=0$: 
$$t =\bigl\{\phi(\bb X):\ \phi(\bb X)\in L,\ \phi(\bb X)
{\rm \ quantifier\ free},\ C\models\phi(\bb A) \bigr\}.$$
for $n=m+1$:
$$t = \bigl\{\th^m(C;\bb A\ee B):\ B\in {\cal P}(C)\}.$$ 
We may regard $\th^n(C;\bb A)$ as the set of $\vf(\bb X)$ that are boolean 
combinations of monadic formulas of quantifier depth $\le n$ \st 
$C\mo \vf(\bb A)$.    \vv
\pd \dd 3.2. $\T_{n,l}$ is the set of all formally possible 
$\th^n(C;\bbp)$ where $C$ is a chain and $\lg(\bbp)=l$.  
$T_{n,l}$ is $|{\cal T}_{n,l}|$.     \vv 
\pt Fact 3.3. (A) For every formula $\psi(\bb X)\in L$ there is an $n$  
such that from $\th^n(C;\bb A)$ we can effectively decide whether 
$C \models \psi(\bb X)$.  If $n$ is minimal with this property we will 
write \/ \underbar{$\dpp(\psi)=n$}.

\\(B) If\/ $m\ge n$ then $\th^n(C;\bb A)$ can be effectively computed 
from $\th^m(C;\bb A)$.  

\\(C) For every  $t\in\T_{n,l}$ there is a monadic 
formula $\psi_t(\bb X)$ with $\dpp(\psi)=n$ \st for every 
$\bb A\in\,^l{\cal P}(C)$, \ \ $C\mo\psi_t(\bb A)\ \iff \ \th^n(C;\bb A)=t$.

\\(D) Each $\th^n(C;\bb A)$ is hereditarily finite, and we can effectively 
compute the set $T_{n,l}$ of formally possible $\th^n(C;\bb A)$. \vv
\pd Proof. Easy. 
\vv \qed
\pd \dd 3.4. If\/ $C,D$ are chains then $C+D$ \ is any chain that can be split
into an initial segment isomorphic to $C$ and a final segment isomorphic to 
$D$.  

\\If\/ $\langle C_i:i<\al\rangle$ is a sequence of chains then 
$\sum_{i<\al}C_i$ \ is any chain $D$ that is the concatenation of segments  
$D_i$, such that each $D_i$ \ is isomorphic to $C_i$. \vv
\pt Theorem 3.5 (composition theorem for linear orders). 

\\(1) If \ $\lg(\bb A)=\lg(\bb B)=\lg(\bb A')=\lg(\bb B')=l$, and
$$\th^m(C;\bb A) = \th^m(C';\bb A') \ \ {\rm and}\ \ 
\th^m(D;\bb B) = \th^m(D';\bb B')$$
then 
$$\th^m(C+D;A_0\cup B_0,\ldots,A_{l-1}\cup B_{l-1}) = 
\th^m(C'+D';A'_0\cup B'_0,\ldots,A'_{l-1}\cup B'_{l-1}).$$ 
(2) If for $i<\al$, \/ $\th^m(C_i;\bb {A_i}) = \th^m(D_i;\bb {B_i})$ \/ 
where \/ $\bb A_i=\lan A_0^i,\ldots,A_{l-1}^i\ran$, 
$\bb B_i=\lan B_0^i,\ldots,B_{l-1}^i\ran$ \ then
$$\th^m\Big(\sum_{i<\al}C_i;\,\cup_{i<\al}A_0^i,\ldots,
\cup_{i<\al}A_{l-1}^i\Big)=\th^m\Big(\sum_{i<\al}D_i;\,\cup_{i<\al}B_0^i,
\ldots,\cup_{i<\al}B_{l-1}^i\Big)$$   \vv
\pd Proof. By \sh\ Theorem 2.4 (where a more general theorem is proved), 
or directly by induction on $m$. 
\vv \qed
\pd \dd 3.6.  (1) \ $t_1+t_2=t_3$ means: 

\\for some $m,l<\om$, \ $t_1,t_2,t_3\in \T_{m,l}$ \/ and if 
$$t_1=\th^m(C;A_0,\ldots,A_{l-1}) \ \ {\rm and} \ \ 
t_2=\th^m(D;B_0,\ldots,B_{l-1})$$ 
then  
$$t_3=\th^m(C+D;A_0\cup B_0,\ldots,A_{l-1}\cup B_{l-1}).$$
By the previous theorem, the choice of\/ $C$ and $D$ is immaterial. 

\\(2) \ $\sum_{i<\al}\th^m(C_i;\bb {A_i})$ \ \ is  \ 
$\th^m(\sum_{i<\al}C_i;\ \cup_{i<\al}A_0^i,\ldots,\cup_{i<\al} A_{l-1}^i)$. 
\pd Notation 3.7. 

\\(1) $\th^n(C;\bb P,\bb Q)$ is $\th^n(C;\bb P\ee\bb Q)$.

\\(2) If\/ $D$ is a subchain of\/ $C$ and $X_1,\ldots,X_{l-1}$ \ are subsets 
of $C$ then $\th^m(D;X_0,\ldots,X_{l-1})$ abbreviates 
$\th^m(D;X_0\cap D,\ldots,X_{l-1}\cap D)$.  

\\(3) For $C$ a chain, $a<b\in C$ and $\bb P\sb C$ we denote by 
$\th^n(C;\bb P)\res_{[a,b)}$ the theory \par
$\th^n([a,b);\bb P\cap[a,b))$.  

\\(4) We will use abbreviations as $\bb P\cup\bb Q$ for 
$\lan P_0\cup Q_0,\ldots\ldots\ran$ and $\cup_i\bb P_i$ for
$\lan \cup_i P_0^i,\ldots\ldots\ran$ (of course we assume that all the 
involved sequences have the same length).  

\\(5) We shall not always distinguish between $\th^n(C;\bbp,\em)$ and
$\th^n(C;\bbp)$. \vv
\pt Theorem 3.8. For every $n,l<\om$ there is $m=m(n,l)<\om$, effectively 
computable from $n$ and $l$, \st whenever $I$ is a chain, for $i\in I\,$ $C_i$ 
is a chain, $\bb Q_i\sb C_i$ and $\lg(\bb(Q_i)=l$,  

\\if $(C;\bb Q)=\sum_{i\in I}(C_i;\bb Q_i):=
(\sum_{i\in I}C_i;\cup_{i\in I}\bb Q_i)$

\\and if for $t\in\T_{n,l}\,$ $P_t:=\{i\in I:\th^n(C_i;\bb Q_i)=t\}$ and 
$\bbp:=\lan P_t:t\in\T_{n,l}\ran$

\\then from $\th^m(I;\bbp)$ we can effectively compute $\th^n(C;\bb Q)$ \vvv
\pd Proof. By \sh\ Theorem 2.4. 
\vv \qed
\pd \dd 3.9. 

\\(1) Let $T_0$, $T_1$ be disjoint trees with $\et_0=root(T_0)$.
Define a tree $T$ to be the ordered sum of \/ $T_0$ and $T_1$ by:

\\$T=T_0\bpl T_1$ iff \/ $T=T_0\cup T_1$ where the 
partial order on $T$, $\seg_T$, \/ is induced by the partial orders of \/ 
$T_0$ and $T_1$ and the (only) additional rule:  
$$\si\in T_1\imp\et_0\seg\si.$$
\\(2) If $T_0$ doesn't have a root then $\seg_T$ is the disjoint union 
$\seg_{T_0}\cup\seg_{T_1}$ (So $[\ta\in T_0\ \&\ \si\in T_1]\imp\ta\bot\si$). 

\\(3) When $I$ is a chain and $T_i$ are pairwise disjoint trees for $i\in I$ 
we define $T=\bpl_{i\in I}T_i$ by $T=\cup_{i\in I}T_i$ with similar rules on 
$\seg=\seg_T$ namely 

$$\si,\ta\in T_i \ \imp \ [\si\seg\ta\iff\si\seg_{T_i}\ta]$$
$$[\si=root(T_i),\, i<_Ij,\, \ta\in T_j] \ \imp \ \si\seg\ta$$ 
$$[\si\in T_i,\, \si\ne root(T_i),\, i\ne j,\, \ta\in T_j] \ \imp \ 
\si\bot\ta$$  \vv
\pt Theorem 3.10 (composition theorem along a complete branch).

\\For every $n<\om$ there is an $m=m(n)<\om$, effectively 
computable from $n$, \st if \/ $I$ is a chain and $T_i$ are trees for $i\in I$ 
then $\lan\th^m(T_i):i\in I\ran$ \/ and \/ $\th^m(\lan \et_i:i\in I\ran)$ \/ 
(which is a theory of a chain) determine $\th^n(\bpl_{i\in I}T_i)$. \vv
\pd Proof. See theorem 3.14.
\vv \qed
Given a tree $T$, we would like to represent it as a sum of subtrees, ordered  
by a branch $B\sb T$. Sometimes however we may have to use a chain $\B$ 
that embeds $B$. 

\pd \dd 3.11. Let $T$ be a tree $T$, $B\sb T$ a branch $\nu\in T$, 
$\et\in B$ and $X\sb B$ be an initial segment without a last element. 

\\(a) $\nu$ {\sl cuts} $B$ at $\et$ if \/ $\et\seg\nu$ and for every 
$\ta\in B$, if $\neg\ta\seg\et$ then $\neg\ta\seg\nu$,
(In particular, $\et$ cuts $B$ at $\et$). $\nu$ cuts $B$ at $\{\et\}$ has 
the same meaning.

\\(b) $\nu$ {\sl cuts} $B$ at $X$ if \/ $\et\seg\nu$ \/ for every $\et\in X$ 
and $\neg\ta\seg\nu$ \/ for every $\ta\in B\sm X$.

\\(c) $\B^+\sb\P(B)$ is defined by $X\in \B^+$ iff \/ $\big[X=\{\et\}$ \/ for 
some $\et\in B\big]$ \ or \ $\big[X\sb B$ is an initial segment without a last 
element and there is $\nu\in T\sm B$ that cuts $B$ at $X\big]$.

\\(d) Define a linear order $\le\,=\,\le_{\B^+}$ on $\B^+$ by $X_0\le X_1$ 
iff $\big[X_0=\{\et_0\},\,X_1=\{\et_1\}$ and $\et_0\seg\et_1\big]$ \ or \ 
$\big[X_0\sb X_1\big]$. 

\\Note that the statements $X\in\B^+$ and $X_0\le_{\B^+}X_1$ are expressible 
by monadic formulas $\psi_{\in}(X,B)$ and $\psi_{\le}(X_0,X_1,B)$.

\\(e) For $X\in\B^+$ define $T_X:=\big\{\nu\in T:\nu$ cuts $B$ at 
$X\big\}$.   \vv
\\

\\Now $\B^+$ has the disatvantage of not being a subset of $T$ and (at the 
small cost of adding a new parameter) we shall replace the chain 
$(\B^+,<_{\B^+})$ by a chain $(\B,<_\B)$ where $\B\sb T$.

\pd \dd 3.12. $\B\sb T$ is obtained by replacing every $X\in \B^+$ by an 
element $\et_{\rm x}\in T$ in the following way: \ if $X=\{\et\}$ then 
$\et_{\rm x}=\et$ and if $X\sb B$ is an initial segment then $\et_{\rm x}$ is 
a favourite element from $T_X$. $\le_\B$ is defined by 
$\et_{{\rm x}_1}\le_\B\et_{{\rm x}_2}\iff 
X_1\le_{\B^+}X_2$ and $B^c\sb T$ will be 
$\B\sm\{\et_X:X=\{\nu\},\,\nu\in B\}$, (so $(\B\sm B^c,\le_\B)\cong(B,\seg)$). 
For $\et\in B$ let $T_\et$ be $T_{\{\et\}}$ as defined in (e) above, and 
for $\et=\et_{\rm x}\in B^c$ let $T_\et=T_X$ as above (in this case $T_\et$ is 
$\{\nu\in T:\nu\sim^0_B\et\}$ as in definition 2.5). \vvv
\pt Fact 3.13. $\le_\B$ is definable from $B$ and $B^c$,  $T_\et$ is definable 
from $\et,B$ and $B^c$ and $T=\bpl_{\et\in\B}T_\et$ in accordance with 
definition 3.9. \vv \qed
\pt Theorem 3.14 (Composition theorems for trees).

\\Assume $T$ is a tree, $B\sb T$ a branch and $\bb Q\sb T$ with 
$\lg(\bb Q)=l$. Let $\B$ and $B^c$ be defined as above, for $\et\in\B$ 
$T_\et$ is defined as above (so $T=\bpl_{\et\in\B}T_\et$) and $S_\et$ is 
$T_\et\sm B$ (so, abusing notations, $T=B\cup\bpl_{\et\in\B}S_\et$). 
\underbar{Then:}

%
%
\\1) \underbar{Composition theorem on a branch:} for every $n<\om$ there is 
\/ $k=k(n,l)<\om$, effectively computable from $n$ and $l$, \st  
$\th^k(\B;B,B^c,\bbp)$ determines $\th^n(T;\bb Q)$ 

\\where for $t\in\T_{n,l}$, $P_t:=\{\et\in\B:\th^n(T_\et;\bb Q\cap T_\et)=t\}$ 
and $\bbp:=\lan P_t:t\in\T_{n,l}\ran$. 

\\2) \underbar{Composition theorem along a branch:} \/ for every $n<\om$ 
there is \/ $k=k(n,l)<\om$, effectively computable from $n$ and $l$, \st 

\\$\th^k(B;\bb Q)$ \/ and \/ $\lan\th^k(S_\et;B,B^c,\bb Q):\et\in\B\ran$ 
determine $\th^n(T;\bb Q)$.            \vv
\pd Proof. By Theorem 1 in \gu $\,\S$2.4.
\vv \qed
\pd \dd 3.15. Additive colouring....
\vv
\pt Theorem 3.16 (Ramsey theorem for additive colourings). ... 
\vv
\pd Proof. By \sh\ Theorem 1.1.
\vv \qed

\sec

%
%
%
%
%
%
\hh 4. Well Orderings of Ordinals \par
%
%
%
\\A chain is {\sl tame} iff it is scattered of Hausdorff degree $<\om$. We 
will define for a tame chain $C$, Log$(C)$ and show later (in proposition 4.8)
that this function is well defined.

\pd \dd 4.1. Let \/ Log:$\{$tame chains$\}\to\om\cup\{\infty\}$ be defined by: 

\\Log$(C)=\infty$ iff there is $\vf(x,y,\bbp)$ that defines a well ordering on 
the elements of $C$ of order type $\ge\om^\om$,

\\Log$(C)=k$ iff there is $\vf(x,y,\bbp)$ that defines a well ordering on 
the elements of $C$ of order type $\al$ with $\om^k\le\al<\om^{k+1}$. \vv
\pt Fact 4.2. A tame chain $C$ has a reconstrutible well ordering i.e. there 
is a formula $\vf(x,y,\bbp)$ ($\bbp\sb C$) that defines a well ordering on the 
elements of $C$ of order type $\al$ and there is a formula $\psi(x,y,\bb Q)$ 
($\bb Q\sb\al$) that defines a linear order $<^*$ on the elements of $\al$ \st 
$(\al,<^*) \cong (C,<)$.
\vv
\pd Proof. By induction on Hdeg$(\al)$, using the proof of Theorem A1 in the 
appendix.
\vv \qed
\pd \dd 4.3. Let $\al,\be$ be ordinals. $\al\to\be$ means the following: 
``there is $\vf(x,y,\bbp)$ that defines a well ordering on the elements of 
$\al$ of order type $\be$''. 
\vv

\pt Claim 4.4. 

\\1) $\al\to\be \ \& \ \be\to\ga \ \imp \ \al\to\ga$.

\\2) $\al\to\ga \ \& \ \ga\ge\al\cdot\om \ \imp \ \al\to\al\cdot\om$.  \vv
\pd Proof. Straightforward.
\vv \qed
\pd Notation. Suppose  $\al\to\be$ holds by $\vf(x,y,\bbp)$. Define   
a bijection $f\colon\al\to\be$ by $f(i)=j$ iff $i$ is the $j$'th element in 
the well order defined by $\vf$.    \vv
\pt Lemma 4.5. For any ordinal $\al$, \ $\al\not\to\al\cdot\om$.   \vv  
\pd Proof. Assume that $\al$ is minimal \st \  $\al\to\al\cdot\om$. It 
follows that:

\\$(i)$ \ \ \  $\al\ge\om$,

\\$(ii)$ \ \  $\al$ is a limit ordinal (by $\al\to\al+1$ and 2.7),

\\$(iii)$ \ for $\be<\al$, \ $\{f(i):i<\be\}$ does not contain a final segment 
of $\al\cdot\om$ (otherwise clearly $\be\to\al\cdot\om$ hence by 2.7 \ 
$\be\to\be\cdot\om$ but $\al$ is minimal).

\\So let $\vf(x,y,\bbp)$ define a well order of $\al$ of order type 
$\al\cdot\om$ and let $Q\sb\al$ be the following subset: \ $x\in Q$ \ iff  
for some $k<\om$, \ $\al\cdot 2k\le f(x)<\al\cdot(2k+1)$. 
\ Let $E$ an equivalence relation on $\al$ defined by $xEy$ \ iff for some 
$l<\om$, \ $f(x)$ and $f(y)$ belong to the segment 
$[\al\cdot l,\al\cdot(l+1))$.  \ 
Clearly there is a monadic formula $\psi(x,y,\bbp,Q)$ that defines $E$ 
moreover, some monadic formula $\theta(X,\bbp,Q)$ expresses the statement  
``\ $\bigvee_{i<\om}\big(X=Q_i\big)$\/'' where $\lan Q_i:i<\om\ran$ are the 
$E$-equivalence classes.

\\Let $n:={\rm max}\big\{\dpp(\vf),\dpp(\psi),\dpp(\theta)\big\}+5$, \ and 

\\$m:=|\big\{\th^n(C;\bb X,Y,Z) \,: \, C {\rm \ a\ chain\ },\ \bb X,Y,Z\sb C, 
\ \lg(\bb X)=\lg(\bbp)\big\}|$.

\\let $\de=\cf(\al)$ and $\{x_i\}_{i<\de}$ be stricly increasing and cofinal 
in $\al$. By \sh Theorem 1.1 applied to the colouring 
$h(i,j)=\th^n(\al;\bbp,Q,x_i,x_j)$ 
we get a cofinal subsequence $\{\be_j\}_{j<\de}$  \st 
$\th^n(\al;\bbp,Q,\be_{j_1},\be_{j_2})$ is constant for $j_1<j_2<\de$. Note 
that it follows  

\\$(\dag)$  the theories $\th^n(\al;\bbp,Q)\res_{[0,{\be_j})}$, \ 
$\th^n(\al;\bbp,Q)\res_{[{\be_j},\al)}$, \ and \
$\th^n(\al;\bbp,Q,\be_{j_1})\res_{[{\be_{j_1}},{\be_{j_2}})}$ 
\ \ are constant for every $j<\de$ and for every $j_1<j_2<\de$.

\\Note that each $E$-equivalence class $Q_i$ is unbounded in $\al$ since if 
some $\be<\al$ contains some $E$-equivalence class $Q_i$ it would easily 
follow that $\be\to\al$ contradicting fact $(iii)$. 

\\Fix some $1<j<\de$ let $x<\be_j$ and let $Q_{i(x)}$ be the $E$-equivalence 
class containing $x$. Since $Q_{i(x)}$ is unbounded in $\al$ there is some 
$j<l<\de$ \st $[\be_j,\be_l)\cap Q_{i(x)}\ne\em$. This statement is 
expressible by \ $\th^n(\al;\bbp,Q,x,\be_j,\be_l)$ which is equal to 

\\$\th^n(\al;\bbp,Q,x,\be_j,\be_l)\res_{[0,{\be_j})}+
\th^n(\al;\bbp,Q,x,\be_j,\be_l)\res_{[{\be_j},{\be_l})}+
\th^n(\al;\bbp,Q,x,\be_j,\be_l)\res_{[{\be_l},\al)}$ \ =

\\$\th^n(\al;\bbp,Q,x,\em,\em)\res_{[0,{\be_j})}+
\th^n(\al;\bbp,Q,\em,\be_j,\em)\res_{[{\be_j},{\be_l})}+
\th^n(\al;\bbp,Q,\em,\em,\be_l)\res_{[{\be_l},\al)}$.

\\By $(\dag)$ we may replace the second theory by 
$\th^n(\al;\bbp,Q,\em,\be_j,\em)\res_{[{\be_j},{\be_{j+1}})}$

\\and the third theory by 
$\th^n(\al;\bbp,Q,\em,\em,\be_{j+1})\res_{[{\be_{j+1}},\al)}$, and conclude:
$$\th^n(\al;\bbp,Q,x,\be_j,\be_l)=\th^n(\al;\bbp,Q,x,\be_j,\be_{j+1})$$
Therefore for every $x<\be_j$, \ $[\be_j,\be_{j+1})\cap Q_{i(x)}\ne\em$. 

Finally, let $j<\de$ be \st the segment $[0,\be_j)$ intersects $m+1$ 
different $E$-equivalence classes, say $Q_{i_0},\ldots,Q_{i_m}$. By the 
previous argument we have $[\be_j,\be_{j+1})\cap Q_{i_l}\ne\em$ \ for every 
$l\le m$. By the choice of $m$ there are different $a,b\in\{i_0,\ldots,i_m\}$ 
\st 

\\$(*)$ \ $\th^n(\al;\bbp,Q,Q_a)\res_{[{\be_j},{\be_{j+1}})} \ = \ 
\th^n(\al;\bbp,Q,Q_b)\res_{[{\be_j},{\be_{j+1}})}$.

\\Let $R\sb\al$ be $\big([0,\be_j)\cap Q_a\big)\cup
\big(([\be_j,\be_{j+1})\cap Q_b\big)\cup\big([\be_{j+1},\al)\cap Q_a\big)$. 

\\Now $\th^n(\al,\bbp,Q,R)=$

\\$\th^n(\al,\bbp,Q,R)\res_{[0,{\be_j})}+
\th^n(\al,\bbp,Q,R)\res_{[{\be_j},{\be_{j+1}})}+
\th^n(\al,\bbp,Q,R)\res_{[{\be_{j+1}},\al)}=$

\\$\th^n(\al,\bbp,Q,Q_a)\res_{[0,{\be_j})}+
\th^n(\al,\bbp,Q,Q_b)\res_{[{\be_j},{\be_{j+1}})}+
\th^n(\al,\bbp,Q,Q_a)\res_{[{\be_{j+1}},\al)}=$  (by $(*)$)

\\$\th^n(\al,\bbp,Q,Q_a)\res_{[0,{\be_j})}+
\th^n(\al,\bbp,Q,Q_a)\res_{[{\be_j},{\be_{j+1}})}+
\th^n(\al,\bbp,Q,Q_a)\res_{[{\be_{j+1}},\al)}=$  

\\$\th^n(\al,\bbp,Q,Q_a)$.

\\But $Q_a$ is an $E$-equivalence class while $R$ is not. Since 
$\th^n(\al,\bbp,Q,Z)$ computes the statement 
``\/$Z$ is $E$-equivalence class\/'' we get a contradiction from  
$\th^n(\al,\bbp,Q,R)=\th^n(\al,\bbp,Q,Q_a)$.     
\vvv \qed
\pt Claim 4.6. If $\al\to\be$ and $\be<\al$ then 
$(\ex\ga_1,\ga_2)\big((\ga_1+\ga_2=\al)\, \&\, (\ga_2+\ga_1=\be)\big)$.
\vv
\pd Proof. Let's prove first:

\\Subclaim: $\om+\om\not\to\om$.

\\Proof of the subclaim: Assume that $\vf(x,y,\bbp)$ well orders $\om+\om$ of 
order type $\om$ and that $\dpp(\vf)=n$, $l(\bbp)=l$. Let $x<^*y$ mean 
$(\om+\om,<)\mo\vf(x,y,\bbp)$.

\line{\hfill{$\to$[Insert\ Ramsey\ theorems]}}
\\Let $\{x_i\}_{i<\om}$ be increasing and unbounded in $[0,\om)$ satisfying, 
for $i<j<\om$ and for some $s_0\in\T_{n,l+2}$ and $t_0\in\T_{n,l+2}$
$$\th^n(\om+\om;x_i,\em,\bbp)\res_{[x_i,x_j)}=s_0, \ \ \ 
\th^n(\om+\om;\em,\em,\bbp)\res_{[x_i,x_j)}=t_0,$$
let $\{y_j\}_{j<\om}$ increasing and unbounded in $[\om,\om+\om)$ satisfying,
for $j<k<\om$ and for some $s_1\in\T_{n,l+2}$ and $t_1\in\T_{n,l+2}$
$$\th^n(\om+\om;\em,y_j,\bbp)\res_{[y_j,y_k)}=s_1, \ \ \ 
\th^n(\om+\om;\em,\em,\bbp)\res_{[y_i,y_k)}=t_1.$$
Using Ramsey Theorem (and as $<^*$ is well founded) we may assume that 
$i_1<i_2\imp x_{i_1}<^*x_{i_2}$ and $j_1<j_2\imp y_{j_1}<^*y_{j_2}$.

\\We will show now that for $0<i<\om$ and $0<j<\om$, 
$\th^n(\om+\om;x_i,y_j,\bbp)$ is constant. Indeed, 

\\$t^*:=\th^n(\om+\om;x_i,y_j,\bbp)=$

\\$\th^n(\om+\om;\em,\em,\bbp)\res_{[0,x_0)}+
\th^n(\om+\om;\em,\em,\bbp)\res_{[x_0,x_i)}+$

\\$\th^n(\om+\om;x_i,\em,\bbp)\res_{[x_i,x_{i+1})}+
\th^n(\om+\om;\em,\em,\bbp)\res_{[x_{i+1},\om)}+
\th^n(\om+\om;\em,\em,\bbp)\res_{[\om,y_0)}+$

\\$\th^n(\om+\om;x_i,\em,\bbp)\res_{[y_i,y_j)}+
\th^n(\om+\om;\em,y_j,\bbp)\res_{[y_j,y_{j+1})}+
\th^n(\om+\om;\em,\em,\bbp)\res_{[y_{j+1},\om+\om)}.$

\\Call the sum $t^*=r_1+r_2+\ldots+r_8$. Now $r_1$ is constant, 
$r_2=t_0\cdot i=t_0$ (check that $t_0+t_0=t_0$), $r_3$ is $s_0$, 
$r_4=t_0\cdot\om$ hence is constant, $r_5$ is constant, $r_6=t_1\cdot j=t_1$, 
$r_7=s_1$ and $r_8=t_0\cdot\om$ hence is constant. Therefore $t^*$ does not 
depend on $i$ and $j$. 

\\Now as $\{y_j\}_{j<\om}$ is unbounded with respect to $<^*$, there is some 
$j<\om$ \st $x_1<^*y_j$. This is expressed by  $\th^n(\om+\om;x_1,y_j,\bbp)$
which we have just seen to be independent of $i$ and $j$ hence 
$$(\all\,0<i<\om)(\all\,0<j<\om)\big[(\om+\om,<)\imp\vf(x_i,y_j,\bbp)\big]$$
it follows that $\ot(\om+\om,<^*)\ge\om+1$, a contradiction. 
This proves $\om+\om\not\to\om$.

\\

\\Returning to the proof of the claim, let $\be$ be the minimal ordinal \st 
there exists some $\al>\be$ with $\al\to\be$ but there aren't any 
$\ga_1,\ga_2\le\al$ with $(\ga_1+\ga_2=\al)\, \&\, (\ga_2+\ga_1=\be)$. Call 
such a $\be$ {\sl weird} and let $\al>\be$ the first ordinal witnessing the 
weirdness of $\be$. By transitivity of \/ $\to$ \/ it is easy to see that 
$\be$ is limit. Moreover, $\ga<\be\imp\be\not\to\ga$ hence if 
$\be=\ga_1+\ga_2$ then $\ga_2+\ga_1\ge\be$. It follows that there are two 
possible cases: either $(*)$ \/ $\ga<\be\imp(\ga+\ga<\be)$, hence 
$\ga<\be\imp(\ga\cdot\om\le\be)$ and $\ga<\be\imp(\ot([\ga,\be))=\be)$, 
or $(**)$ \/ $\be=\ga+\ga$.

\\

\\First case: $(*)$ holds i.e. $\ga<\be\imp(\ga+\ga<\be)$. Let $\al=\be+\ga$ 
what can $\ga$ be? If $\ga<\be$ then by $(*)$ $\ga+\be=\be$ and $\al$ does not 
witness the weirdness of $\be$, so $\al\ge\be+\be$. 

\\Let $\vf(x,y,\bbp)$ well order $\al$ of order type $\be$ with $\dpp(\vf)=n$ 
and $l(\bbp)=l$. As above $x<^*y$ means $(\al,<)\mo\vf(x,y,\bbp)$ and finally 
let $\de=\cf(\be)$.  

\\Now $\ot(\al,<^*)=\be$ but what is $\ot([0,\be),<^*\res_{[0,\be)})$? 
Clearly, as $\th^n(\al,\bbp)=
\th^n(\al,\bbp)\res_{[0,\be)}+\th^n(\al,\bbp)\res_{[\be,\al)}$ we have 
$\be\to \ot([0,\be),<^*\res_{[0,\be)})$ hence 
$\be=\ot([0,\be),<^*\res_{[0,\be)})$ (otherwise, by $(*)$,  
$\ot([0,\be),<^*\res_{[0,\be)})$ is weird and $<\be$).
Similarly we can show that $\ot([\be,\be+\be),<^*\res_{[\be,\be+\be)})=\be$. 

\line{\hfill{$\to$[Insert\ Ramsey\ theorems]}}
\\Now proceed as before: choose $\{x_i\}_{i<\de}\sb[0,\be)$ and 
$\{y_j\}_{j<\de}\sb[\be,\be+\be)$ that are homogeneous unbounded and $<^*$ 
unbounded and use them to show that $\ot(\al,<^*)\ge\be+1$.

\\

\\Second case: $(**)$ holds i.e. $\be=\ga+\ga$. 

\\Call $\ep$ {\sl quite weird} if for some $k<\om$ \/ $\ep\cdot k$ is weird. 
Let $\ep\le\ga$ be the first quite weird ordinal. Let $k_1$ be the first \st
$\ep\cdot k_1$ is weird. Look at $\ga$: if $\ga=\ga_1+\ga_2$ and 
$\ga_2+\ga_1<\ga$ we would have $\al\to\be=\ga+\ga\to\ga+\ga_2+\ga_1<\be$ 
and a contradiction. Hence either $\ga_1<\ga\imp(\ga_1+\ga_1<\ga)$ and in this 
case $\ga=\ep$ or $\ga=\ga_1+\ga_1$. Repeat the same argument to get 
$\ga_1=\ep$ or $\ga_1=\ga_2+\ga_2$. After finitely many steps we are bound to   
get $\be=\ep\cdot2k$ where $2k=k_1$ and $\ep_1<\ep\imp\ep_1\cdot\om\le\ep$ 
and of course $\ep_1<\ep\imp\ep\not\to\ep_1$.

\\Let $\vf(x,y,\bbp)$ and $<^*$ be as usual and $\de:=\cf(\be)=\cf(\ep)$.  
Let $\al=\be+\ep^*$ if $\ep^*<\ep$ then $\ep^*+\be=\be$ and $\al$ doesn't 
witness weirdness, therefore $\ep^*\ge\ep$.

\\Proceed as before: choose $\{x^0_i\}_{i<\de},\{x^1_i\}_{i<\de},\ldots, 
\{x^k_i\}_{i<\de}$ with $\{x^l_i\}_{i<\de}\sb[\ep\cdot l,\ep(l+1))$, 
homogeneous, unbounded and $<^*$ incresing. 

\\By the composition theorem it will follow that 
$\ot([\ep\cdot l,\ep(l+1)),<^*)\ge\ep$ and by homogeneity we will have, for 
$0<i,j<\om$ and $l\le k$, $x^l_i<^*x^{l+1}_j$.
It follows that $\ot(\al,<^*)\ge(\ep\cdot k)+1=\be+1$ and a contradiction.
\vv \qed
\pt Theorem 4.7. Well ordering of ordinals are obtained only by the following 
process: 

\\let $\lan P_0,P_1,\ldots,P_{n-1}\ran$ be a partition of $\al$ and 
$$i<^*j \iff \big[(\ex k<n)(i\in P_k\,\&\,j\in P_k\,\&\,i<j)\big]\vee
\big[i\in P_{k_1}\,\&\,j\in P_{k_2}\,\&\,k_1<k_2\big].$$
\vv \qed
\pt Proposition 4.8. Log$(C)$ is well defined. \vv
\pd Proof. Let $(C,<^*)$ be a scattered chain and let $(\al,<)$ and  
$(\be,<)$ be results of a definable well orderings of\/ $(C,<^*)$ where in 
addition (by 4.2) there is $\psi(x,y,\bb Q)$ that defines $C$ in $\al$.  
So $\al\to\be$ and by 4.5 and 4.6 $\al<\om^\om \iff \be<\om^\om$ and    
$\al\in[\om^k,\om^{k+1}) \iff \be\in[\om^k,\om^{k+1})$.
\vv \qed
\sec 
%
%
%
%
%
%
\hh 5. $(\om^\om,<)$ and longer chains \par
%
%
%
\\The following lemma is a part of Theorem 3.5(B) in \sh:

\pt Lemma 5.1. Let $I$ be a well ordered chain of order type $\ge\om^k$.  
Let $f\colon I^2\to\{t_0,t_1\ldots,t_{l-1}\}$ be an additive colouring and  
assume that for $\al<\be\in I$, $f(\al,\be)$ depends only on the order type in 
$I$ of the segment $[\al,\be)$.  

\\Then there is \/ $i<l$ \st for some $p\le l$, for every $r\ge p$, if 
$\ot\big([\al,\be)\big)=\om^r$ then $f(\al,\be)=t_i$. Moreover, $t_i+t_i=t_i$.
\vvv
\pd Proof. To avoid triviality assume $k>l$. For $\al<\be$ in $I$ with 
$\ot\big([\al,\be)\big)=\de$, denote $f(\al,\be)$ by $t(\de)$  (makes sense 
by the assumptions). 

\\By the pigeon-hole principle there are $1\le p\le l$, $s>p$ and some $t_i$ 
with $t(\om^p)=t(\om^s)=t_i$. Now $\om^{p+2}=\sum_{i<\om}(\om^{p+1}+\om^p)$ and 
by the additivity of $f$:
$$t\big(\om^{p+2}\big)=t\big(\sum_{i<\om}(\om^{p+1}+\om^p)\big)=
\sum_{i<\om}t\big(\om^{p+1}+\om^p\big)=
\sum_{i<\om}\big(t(\om^{p+1})+t(\om^p)\big)=
\sum_{i<\om}\big(t(\om^{p+1})+t(\om^s)\big)=$$
$$\sum_{i<\om}t\big(\om^{p+1}+\om^s\big)=
\sum_{i<\om}t\big(\om^s\big)=\sum_{i<\om}t\big(\om^p\big)=
t\big(\sum_{i<\om}\om^p\big)=t\big(\om^{p+1}\big).$$
Hence 
$$t\big(\om^{p+2}\big)=t\big(\om^{p+1}\big).$$
\\

\\Using this and as $\om^{p+3}=\sum_{i<\om}(\om^{p+2}+\om^{p+1})$ we have 
$$t\big(\om^{p+3}\big)=t\big(\sum_{i<\om}(\om^{p+2}+\om^{p+1})\big)=
\sum_{i<\om}t\big(\om^{p+2}+\om^{p+1}\big)=
\sum_{i<\om}\big(t(\om^{p+2})+t(\om^{p+1})\big)=$$
$$\sum_{i<\om}\big(t(\om^{p+1})+t(\om^{p+1})\big)=
\sum_{i<\om}t\big(\om^{p+1}\big)=
t\big(\sum_{i<\om}\om^{p+1}\big)=t\big(\om^{p+2}\big).$$
Hence 
$$t\big(\om^{p+3}\big)=t\big(\om^{p+2}\big).$$
\\

\\So for every $j>0$, \  $t\big(\om^{p+1}\big)=t\big(\om^{p+j}\big)$ and in 
particular \ $t\big(\om^{p+1}\big)=t\big(\om^s\big)=t\big(\om^p\big)=t_i$.

\\This proves the first part of the lemma. As for the moreover clause, 
since $\om^{p+1}=\om^p+\om^{p+1}$ \ we have
$$t_i=t\big(\om^{p+1}\big)=t\big(\om^p+\om^{p+1}\big)=
t\big(\om^p\big)+t\big(\om^{p+1}\big)=t_i+t_i.$$
\vv \qed
\pt Proposition 5.2. The formula $\vf(X,Y)$ saying ``if\/ $Y$ is without a 
last element then $X\sb Y$ is an $\om$-sequence unbounded in $Y$  
(and if not then $X=\em$)'' can not be uniformized in $(\om^\om,<)$. 

\\Moreover, if \/ $\psi_m(X,Y,\bb P_m)$ uniformizes $\vf$ on $\om^m$ then one of 
the sets $\{\dpp(\psi_m):m<\om\}$ or $\{\lg(\bb P_m):m<\om\}$ is unbounded. \vvv
\pd Proof. Suppose the second statement fails, then:

\\$(\dag)$ there is a formula $\psi(X,Y,\bb Z)$ \st for an unbounded set 
$I\sb\om$, for every $m\in I$ there is $\bb P_m\sb\om^m$ \st 
$\psi(X,Y,\bb P_m)$ uniformizes $\vf$ on $\om^m$.

\\Let $\bb P_m=\bbp$ let $n=\dpp(\psi)+1$ and 
$M:=|\big\{\th^n(C;X,Y,\bb Z) \,: \, C {\rm \ a\ chain\ },\ X,Y,\bb Z\sb C, 
\ \lg(\bb Z)=\lg(\bbp)\big\}|$.

\\Let $m\in I$ be large enough ($m>2M+3$ will do), and let's show that $\psi$ 
doesn't work for $\om^m$ and a subset $Y_1$ that will be defined now.

\\If $\al<\om^m$ then 
$\al=\om^{m-1}k_{m-1}+\om^{m-2}k_{m-2}+\ldots+\om k_1+k_0$. Let $k(\al):=
{\rm min}\{i:k_i\ne 0\}$ and let $A_k:=\{\al<\om^m:k(\al)=k\}$. Note that 
$\ot(A_k)=\om^{m-k}$.

\\For $k\in\{1,2,\ldots,m-1\}$ we will choose $Y_k\sb A_k$ with $\ot(Y_k)=
\ot(A_k)=\om^{m-k}$ \st for $\al<\be$ in $Y_k$:
$$(*) \ \ \ \th^n(\om^m;\bbp,Y_k)\res_{[\al,\be)}\ {\rm depends\ only\ on\ } 
\ot\big([\al,\be)\cap Y_k\big)$$
we will start with $k=m-1$ and proceed by inverse induction:

\\Let $A_{m-1}=\lan\al_j:j<\om\ran$.
Let for $l<p<\om$, \ $h(l,p):=\th^n(\om^m;\bbp,\al_l)\res_{[{\al_l},{\al_p})}$.
Let $J\sb\om$ be homogeneous with respect to this colouring namely, for some 
fixed theory $t_{m-1}$, for every $l<p$ in $J$,  
$$\th^n(\om^m;\bbp,\al_l)\res_{[{\al_l},{\al_p})}=t_{m-1}.$$
By the composition theorem, for every $l<p$ in $J$,
$$\th^n(\om^m;\bbp,Y_{m-1})\res_{[{\al_l},{\al_p})}=
t_{m-1}\cdot|Y_{m-1}\cap[\al_l,\al_p)|$$ 
and this proves $(*)$ for $Y_{m-1}$.

\\

\\Rename $Y_{m-1}$ by $\lan\al_i:i<\om\ran$. 
In each segment $[\al_i,\al_{i+1})$ choose 
$\lan\be_l^i:0<l<\om\ran\sb A_{m-2}$  increasing and cofinal \st for every 
$l<p<\om$ the theory $\th^n(\om^m;\bbp,\be_l^i)\res_{[{\be_l^i},{\be_p^i})}$ 
is constant.

\\Returning to $Y_{m-1}$, for $i<j<\om$ let 
$$h_1(i,j):=\big\lan\th^n(\om^m;\bbp)\res_{[{\al_i},{\be_1^{j-1}})}\, ,
\, \th^n(\om^m;\bbp,\be_1^{j-1})\res_{[{\be_1^{j-1}},{\be_2^{j-1}})}\big\ran$$
w.l.o.g. (by thinning out and re-renaming and noting that we don't harm 
$(*)$) $Y_{m-1}$ is homogeneous with respect to this colouring. 

\\Hence, for some theories $t^*$ and $t_{m-2}$, for every $i<j<\om$ we have 
$$h_1(i,j)=\lan t^*,t_{m-2}\ran$$
Let $Y_{m-2}:=\lan\be_l^i:0<l<\om,i<\om\ran$, clearly  
$\ot(Y_{m-2})=\om^2$. \ Let's check $(*)$ for $Y_{m-2}$:

\\Firstly, note that for $l<p<\om$, 
$$\th^n(\om^m;\bbp,Y_{m-2})\res_{[{\be_l^i},{\be_p^i})}=t_{m-2}\cdot(p-l).$$
Secondly, for $i<j<\om$ \ \ 
$\th^n(\om^m;\bbp,Y_{m-2})\res_{[{\be_l^i},{\be_p^j})}=$
$$\th^n(\om^m;\bbp,Y_{m-2})\res_{[{\be_l^i},{\al_{i+1}})}+
\th^n(\om^m;\bbp,Y_{m-2})\res_{[{\al_{i+1}},{\al_{i+2}})}+\ldots \ \ \ \ $$
$$\ \ \ \ +\th^n(\om^m;\bbp,Y_{m-2})\res_{[{\al_{j-1}},{\al_j})}+
\th^n(\om^m;\bbp,Y_{m-2})\res_{[{\al_j},{\be_p^j})}$$
where the first theory is equal to $t_{m-2}\cdot\om$, the last theory is 
$t^*+t_{m-2}\cdot(p-l)$, and the middle theories are $t^*+t_{m-2}\cdot\om$. 
These observations prove $(*)$ for $Y_{m-2}$.

\\

\\For defining $Y_{m-3}$ let's restrict ourselves to a segment 
$[\al_i,\al_{i+1})$ where $\al_i,\al_{i+1}\in Y_{m-1}$. In this segment we 
have defined $\lan\be^i_l:0<l<\om\ran\sb Y_{m-2}$. Now choose in each 
$[\be^i_l,\be^i_{l+1})$ an increasing cofinal sequence 
$\lan\ga^{i,l}_j:0<j<\om\ran$ \st for $j<p<\om$,  \ 
$\th^n(\om^m;\bbp,\ga^{i,l}_j)\res_{[{\ga^{i,l}_j},{\ga^{i,l}_p})}$ is 
constant.

\\For \/ $0<l<p<\om$ let \/ 
$$h_1^i(l,p):=
\big\lan\th^n(\om^m;\bbp)
\res_{[{\be^i_l},{\ga^{i,p-1}_1})}, \,  
\th^n(\om^m;\bbp,\ga^{i,p-1}_1)
\res_{[{\ga^{i,p-1}_1},{\ga^{i,p-1}_2})}\big\ran$$ 
and again w.l.o.g we may assume that $\lan\be^i_l:0<l<\om\ran$ is homogeneous 
with respect to $h_1^i$.

\\Next, for $i<j<\om$ define 
$$h_2(i,j):=
\big\lan\th^n(\om^m;\bbp)
\res_{[{\al_i},{\ga^{j-1,1}_1})}, \,  
\th^n(\om^m;\bbp,\ga^{j-1,1}_1)
\res_{[{\ga^{j-1,1}_1},{\ga^{j-1,1}_2})}\big\ran$$ 
by thinning out and renaming we may assume that $Y_{m-1}$ is homogeneous with 
respect to $h_2$, now $Y_{m-2}$ is also thinned out but each new  
$\lan\be^i_l:0<l<\om\ran$ which is some old $\lan\be^{i^*}_l:0<l<\om\ran$ 
is still homogeneous.

\\As a result we will have, for some theories $t^{**},t^{***},t_{m-3}$:
$$(\all i<j<\om)(\all 0<l<p<\om)\big[h_1^i(l,p)=\lan t^{**},t_{m-3}\ran \ \&
\ h_2(i,j)=\lan t^{***},t_{m-3}\ran\big].$$
Let $Y_{m-3}:=\{\ga^{i,l}_j:i<\om,0<l<\om,0<j<\om\}$, as before $(*)$ holds 
by noting that if for example $i_1<i_2<\om$ and $1<l_2$ then
$$\th^n(\om^m;\bbp,\ga^{{i_1},{l_1}}_{j_1})
\res_{[{\ga^{{i_1},{l_1}}_{j_1}},{\ga^{{i_2},{l_2}}_{j_2}})} \ = \    
t_{m-3}\cdot\om+(t^{**}+t_{m-3}\cdot\om)\cdot\om+
[t^{***}+(t^{**}+t_{m-3}\cdot\om)\cdot\om]\cdot(i_2-i_1-1)+$$
$$t^{***}+t_{m-3}\cdot\om+(t^{**}+t_{m-3}\cdot\om)(l_2-1)+t^{**}+
t_{m-3}\cdot(j_2-1)$$
and similarly for the other possibilities.
 
\\

\\$Y_{m-4},Y_{m-5},\ldots,Y_1$ are defined by using the same prescription i.e. 
$Y_{m-l}$ is defined by taking a homogenous sequence between two successive 
elements of $Y_{m-l-1}$ then homogenous sequences between two successive
elements of $Y_{m-l-2}$ by using colouring of the form $h_1,h_2,\ldots$. 
The thinning out and w.l.o.g's for already defined $Y_{m-k}$'s are not 
necessary but they ease notations considerably.

\\

\\We will show now that $\psi$ doesn't choose an unbounded $\om$-sequence in 
$Y_1$  that is, for every $\om$-sequence $X\sb Y_1$ there is an $\om$-sequence
$X'\sb Y_1$ \st $\th^{n-1}(\om^m;\bbp,Y_1,X)=\th^n(\om^m;\bbp,Y_1,X')$.

\\By $(*)$, for $\al<\be$ in $Y_1$ the additive colouring $f(\al,\be):=
\th^n(\om^m;\bbp,Y_1)\res_{[\al,\be)}$ depends only on 
$\ot\big([\al,\be)\cap Y_1\big)$ hence we can apply lemma 5.1 and conclude 
that for some $p\le m/2$, for every $r\ge p$, \ 
$\th^n(\om^m;\bbp,Y_1)\res_{[\al,\be)}$ is equal to some fixed theory $t$ 
whenever $\ot\big([\al,\be)\cap Y_1\big)=\om^r$. (Remember that $f$ has at 
most $M$ possibilities and that $m>2M$). Moreover, we know that $t+t=t$.

\\Assume now that for some $X\sb Y_1$, \/ $\psi(X,Y_1,\bbp)$ holds, so $X$ is  
a cofinal $\om$-sequence. Let $X=\{\de_i:i<\om\}$. As $\ot(Y_1)=\om^{m-1}$ 
\/ for unboundedly many $i$'s we have 
$\ot\big([\de_i,\de_{i+1})\cap Y_1\big)\ge\om^{m-2}>\om^p$. 

\\Let $\be_i:=\ot\big([\de_i,\de_{i+1})\cap Y_1\big)$ and denote by $t(\ep)$ 
the theory $\th^n(\om^m;\bbp,Y_1)\res_{[\al,\be)}$ when 
$\ot\big([\al,\be)\cap Y_1\big)=\ep$ (by $(*)$ it doesn't matter which $\al$ 
and $\be$ we use).

\\We are interested in $\th^{n-1}(\om^m;\bbp,Y_1,X)$ which is 

\\$\th^{n-1}(\om^m;\bbp,Y_1,\em)\res_{[0,{\de_0})}+
\th^{n-1}(\om^m;\bbp,Y_1,\de_0)\res_{[{\de_0},{\de_1})}+
\th^{n-1}(\om^m;\bbp,Y_1,\de_1)\res_{[{\de_1},{\de_2})}+\ldots$.

\\As $\de_i$ is the first element in $[\de_i,\de_{i+1})\cap Y_1$, \ 
$\th^{n-1}(\om^m;\bbp,Y_1,\de_i)\res_{[{\de_i},{\de_{i+1}})}$ is 
determined by

\\$\th^n(\om^m;\bbp,Y_1)\res_{[{\de_i},{\de_{i+1}})}=t(\be_i)$ and abusing 
notations we will say
$$(**) \ \ \th^{n-1}(\om^m;\bbp,Y_1,X)\ \simeq \ 
t(\de_0)+\sum_{i<\om}t(\be_i).$$
Let $i<\om$ be \st $\be_i\ge\om^{m-2}$ and let $j>i$ be the first with  
$\be_j\ge\om^{m-2}$. 

\\First case: $i=j+1$. 

\\Let $\be_i=\ot\big([\de_i,\de_{i+1})\cap Y_1\big)=\om^{m-2}\cdot k_1+\ep_1$ 
\ and \ $\be_{i+1}=\ot\big([\de_{i+1},\de_{i+2})\cap Y_1\big)=
\om^{m-2}\cdot k_2+\ep_2$ \/ where $k_1,k_2\ge 1$ and $\ep_1,\ep_2<\om^{m-2}$. 

\\Define $\ga:=$ \/ the $\om^{m-2}\cdot k_1+\om^{m-3}+\ep_1$'th 
successor of $\de_i$ in $Y_1$. 
So $\de_{i+1}<\ga<\de_{i+2}$ \ but 
$\ot\big([\de_{i+1},\de_{i+2})\cap Y_1\big)=\be_{i+1}$ hence 
$$\th^n(\om^m;\bbp,Y_1)\res_{[\ga,{\de_{i+2}})}=
\th^n(\om^m;\bbp,Y_1)\res_{[{\de_{i+1}},{\de_{i+2}})}=t(\be_{i+1})$$
hence
$$\th^{n-1}(\om^m;\bbp,Y_1,\ga)\res_{[\ga,{\de_{i+2}})}=
\th^{n-1}(\om^m;\bbp,Y_1,\de_{i+1})\res_{[{\de_{i+1}},{\de_{i+2}})}.$$
On the other hand, 
$$\th^n(\om^m;\bbp,Y_1)\res_{[{\de_i},\ga)}=t(\om^{m-2}\cdot k_1)+
t(\om^{m-3})+t(\ep_1)$$
but $m-3\ge p$ hence $t(\om^{m-3})=t(\om^{m-2})=t$ moreover $t+t=t$ and it
follows that 
$$t(\om^{m-2}\cdot k_1)+t(\om^{m-3})=
t(\om^{m-2})\cdot k_1+t(\om^{m-3})=
t(\om^{m-2})\cdot(k_1+1)=
t(\om^{m-2})\cdot(k_1)=t(\om^{m-2}\cdot k_1)$$
hence
$$\th^n(\om^m;\bbp,Y_1)\res_{[{\de_i},\ga)}=t(\om^{m-2}\cdot k_1)+
+t(\ep_1)=
\th^n(\om^m;\bbp,Y_1)\res_{[{\de_i},{\de_{i+1}})}=t(\be_{i+1})$$
hence
$$\th^{n-1}(\om^m;\bbp,Y_1,\de_i)\res_{[{\de_i},\ga)}=
\th^{n-1}(\om^m;\bbp,Y_1,\de_i)\res_{[{\de_i},{\de_{i+1}})}.$$
Now all other relevant theories are left unchanged therefore, letting 
$X':=X\sm\{\de_{i+1}\}\cup\{\ga\}$ we get $X\ne X'$ but 
$$\th^{n-1}(\om^m;\bbp,Y_1,X)=\th^n(\om^m;\bbp,Y_1,X')$$.

\\General case: $j=i+l$.

\\Look at $\de_{i+1},\de_{i+2},\ldots,\de_{i+l-1},\de_{i+l}=\de_j$. We'll 
define $\ga_1,\ga_2,\ldots,\ga_l$ \/ with $\de_{i+k}<\ga_k<\de_{i+k+1}$ 
for $0<k<l$ and $\ga_l=\de_{i+l}=\de_j$. This will be done by `shifting' 
the $\de_{i+k}$'s by $\om^{m-3}$ (remember that \/ $\be_{i+k}<\om^{m-2}$ 
for $0<k<l$).

\\Assume as before that $\be_i=\ot\big([\de_i,\de_{i+1})\cap Y_1\big)=
\om^{m-2}\cdot k_1+\ep_1$ \/ where $k_1\ge 1$ and $\ep_1<\om^{m-2}$. 

\\Define $\ga_1:=$ \/ the $\om^{m-2}\cdot k_1+\om^{m-3}+\ep_1$'th 
successor of $\de_i$ in $Y_1$, 
$\ga_2:=$ \/ the $\be_{i+1}$'th successor of $\ga_1$ in $Y_1$,
$\ga_3:=$ \/ the $\be_{i+2}$'th successor of $\ga_2$ in $Y_1$ and so on,  
$\ga_l$ will clearly be equal to $\de_j$.

\\As before we have for $1<k\le l$, (by preserving the order types)
$$\th^{n-1}(\om^m;\bbp,Y_1,\ga_k)\res_{[{\ga_k},{\ga_{k+1}})}=
\th^{n-1}(\om^m;\bbp,Y_1,\de_{i+k})\res_{[{\de_{i+k}},{\de_{i+k+1}})}.$$
and (using $t+t=t$)
$$\th^{n-1}(\om^m;\bbp,Y_1,\de_i)\res_{[{\de_i},\ga_1)}=
\th^{n-1}(\om^m;\bbp,Y_1,\de_i)\res_{[{\de_i},{\de_{i+1}})}.$$
Letting $X':=X\sm\{\de_{i+1},\de_{i+2},\ldots,\de_{j-1}\}\cup\{\ga_1,\ga_2,
\ldots,\ga_{l-1}\}$ we get $X\ne X'$ but 
$$\th^{n-1}(\om^m;\bbp,Y_1,X)=\th^n(\om^m;\bbp,Y_1,X')$$.

Since $\dpp(\psi)=n-1$, \/ $X$ is not the unique $\om$-sequence chosen by  
$\psi$ from $Y_1$. Therefore, $\psi$ does not uniformize $\vf$ on $\om^m$,  
a contradiction.  

\\[complete, using composition theorem, for $\om^\om$]
\vvv \qed
\pt Theorem 5.3. If $C$ has the uniformization property then Log$(C)<\om$. 
\vv \qed
\sec
%
%
%
%
%
%
\hh 6. Very Tame Trees   \par
%
%
%
\pt Proposition 6.1. If the ordinals $\al$ and $\be$ have the uniformization 
property then so do $\al+\be$ and $\al\be$. \vvv
\pd Proof. $\al+\be$ is similar to $\al+\al=\al\cdot2$ \/ and we leave it to 
the reader. We shall prove that $\al\cdot\be$ has the uniformization property.

\\Let $\vf(X,Y,\bb Q)$ be p.u in $\al\be$ with $\dpp(\vf)=n$ and
$\lg(\bb Q)=l$. Let $\lan t_0,\ldots,t_{a-1}\ran$ be an enumeration of the 
the theories in $\T_{n,l+2}$. For $i<a$ and $X,Y\sb\al\be$ 
define $P_i(X,Y,\bb Q)\sb K:=\{\al\ga:\ga<\be\}$ by

$$P_i(X,Y,\bb Q):=\big\{\al\ga:
\th^n(\al\be;X,Y,\bb Q)\res_{[\al\ga,\al\ga+\al)}=t_i\big\}$$
it follows that, for every $X,Y\sb\al\be$, \/ 
$\bbp=\bbp(X,Y,\bb Q)=\lan P_0(X,Y,\bb Q),\ldots,P_{a-1}(X,Y,\bb Q)\ran$ 
is a partition of $K$ that is definable from 
$X,Y,\bb Q$ and $K$.

\\$\al\cdot\be$=$\sum_{\ga<\be}[\al\ga,\al\ga+\al)$ and by theorem 3.8 
there is $m=m(n,l)$ \st $\th^n(K;\bbp(X,Y,\bb Q))$ determines 
$\th^n(\al\be;X,Y,\bb Q)$. 

\\Let ${\cal R}=\{r_0,\ldots,r_{c-1}\}$ be the set of theories that satisfy, 
for every $X,Y\sb\al\be$:
$$\th^n(K;\bbp(X,Y,\bb Q))\in{\cal R} \imp \al\be\mo\vf(X,Y,\bb Q).$$
Now let $\lan s_0,\ldots,s_{b-1}\ran$ be an enumeration of the 
the theories in $\T_{n+1,l+1}$. For $i<b$ and $Y\sb\al\be$ define 
$R^0_i(Y,\bb Q)\sb K$ by
$$R^0_i(Y,\bb Q):=\big\{\al\ga:
\th^{n+1}(\al\be;Y,\bb Q)\res_{[\al\ga,\al\ga+\al)}=s_i\big\}$$
as before, for every $Y\sb\al\be$, \/ 
$\bb R^0=\bb R^0(Y,\bb Q)=\lan R^0_0(Y,\bb Q),\ldots,R^0_{b-1}(Y,\bb Q)\ran$ 
is a partition of $K$ that is definable from $Y,\bb Q$ and $K$.

\\Now let $\bb R^1=\lan R^1_0,\ldots,R^1_{a-1}\ran$ be any partition of $K$.
We will say that $\bb R^0(Y,\bb Q)$ and $\bb R^1$ are \underbar{coherent} if 

\\(1) $\al\ga\in(R^0_i\cap R^1_j)$ implies that for every chain $C$, 
$B\sb C$ and $\bb D\sb C$ of length $l$:

if $\th^{n+1}(C;B,\bb D)=s_i$ then 
$(\ex A\sb C)\big[\th^n(C;A,B,\bb D)=t_j\big]$,

\\(2) $\th^n(K;\bb R^1)\in{\cal R}$.

\\Since $a,b$ and $c$ are finite, there is a formula $\theta_1(\bb U,\bb W)$ 
(with $\lg(\bb U)=b$ and $\lg(\bb W)=a$) \st for any $\bb R^0,\bb R^1\sb K$,  

\\$K\mo\theta_1(\bb R^0,\bb R^1)$ \/ iff \/ $\bb R^0$ and $\bb R^1$ are 
coherent partitions of $K$.

\\Moreover, as $K\cong\be$ and $\be$ has the uniformization property, there  
exists $\bb S\sb K$ and a formula $\theta_2(\bb U,\bb W,\bb S)$ \st for every 
$\bb R^0\sb K$

\\if $(\ex\bb W)\theta_1(\bb R^0,\bb W)$ then 
$(\ex!\bb W)[\theta_2(\bb R^0,\bb W,\bb S)\,\&\,\theta_1(\bb R^0,\bb W)].$
Let $\theta(\bb U,\bb W,\bb S):=\theta_1\wedge\theta_2$.

\\Now let $Y\sb\al\be$, let $\bb R^0=\bb R^0(Y,\bb Q)$ and suppose that 
$\bb R^0$ and some $\bb R^1$ are coherent partitions of $K$. 
When $\al\ga\in(R^0_i\cap R^1_j)$, we know by the first clause in the 
definition of coherence that  

\\$(\ex X\sb \al\be)
\big[\th^n(\al\be;X,Y,\bb Q)\res_{[\al\ga,\al\ga+\al)}=t_j\big]$. 

\\Now as $[\al\ga,\al\ga+\al)\cong\al$ and $\al$ has the uniformization 
property, there is $\bb T_\ga\sb[\al\ga,\al\ga+\al)$ and a formula 
$\psi^\ga_j(X,Y,\bb T_\ga)$ (of depth $k(n,l)$ that depends only on $n$ and 
$l$) that uniformizes the formula that says 
``$\th^n(\al\be;X,Y,\bb Q)\res_{[\al\ga,\al\ga+\al)}=t_j$''.  

\\It follows that when $\psi^\ga_j(X,Y,\bb T_\ga)$ holds, 
$X\cap[\al\ga,\al\ga+\al)$ is unique.   

\\W.l.o.g all $\bb T_\ga$ have the same length and (by taking prudent 
disjunctions) $\psi^\ga_j(X,Y,\bb T_\ga)=\psi_j(X,Y,\bb T_\ga)$ and let 
$\bb T=\cup_{\ga<\be}\bb T_\ga$ (the union is disjoint). \ We are ready to 
define $U(X,Y,\bb Q,\bb T,\bb S)$ that uniformizes $\vf(X,Y,\bb Q)$:

\\$U(X,Y,\bb Q,\bb T,\bb S)$ says: ``for every partition $\bb R^0$
of $K$ that is equal to [the definable] $\bb R^0(Y,\bb Q)$ every $\bb R^1$ 
that is a [in fact the only] partition that satisfies 
$\theta(\bb R^0,\bb R^1,\bb S)$, {\sl if} $\al\ga\in R^1_j$ and 
$D=[\al\ga,\al\ga+\al)$ [$\al\ga$ and $\al\ga+\al$ are two successive elements 
of $K$] {\sl then} $D\mo\psi_j(X\cap D,Y\cap D,\bb Q\cap D,\bb T\cap D)$''. 

\\Check that $U(X,Y,\bb Q,\bb T,\bb S)$ does the job: clause (1) in the 
definition of coherence and the $\psi_j$'s guarantee that $X$ is unique, 
clause (2) guarantees that $U(X,Y,\bb Q,\bb T,\bb S)\imp \vf(X,Y,\bb Q)$.
\vv \qed
\pt Fact 6.2. Every finite chain has the uniformization property.
\vv \qed
\pt Theorem 6.3. $(\om,<)$ has the uniformization property.  \vvv
\pt Corollary 6.4. An ordinal $\al$ has the uniformization property iff 
$\al<\om^\om$.   \vvv
\pd \dd 6.5. $(T,\seg)$ is {\sl very tame} if 

\\1) $T$ is tame

\\2) $Sup\{{\rm Log}(B): B\sb T,\ B\ {\rm a\ branch}\}<\om$ 
\vv
\pt Lemma 6.6. If $(T,\seg)$ is not very tame then $(T,\seg)$ does'nt have 
the uniformization property. \vv
\pd Proof. If $T$ is not tame then by theorem 2.7 it doesn't have even a
definable choice function. 

\\If $T$ is tame then either there is a a branch $B\sb T$ with Log$(B)=\infty$ 
or it has branches of unbounded Log. By 3.14(3) and 5.2 and using the 
definable well ordering of $T$, there is a formula $\vf(X,Y,Z)$ that can't be 
uniformized. 
\vv \qed
\pt Theorem 6.7. $(T,\seg)$ has the uniformization property iff\/ $(T,\seg)$ 
is very tame.   \vv
\pd Proof. Assume $T$ is $(l^*,n^*,k^*)$ very tame and let $\vf(X,Y,\bb Q)$ 
be p.u in $T$ with $\dpp(\vf)=n$ and $\lg(\bb Q)=l$. 

\\As $T$ is $(n^*,k^*)$ tame it can be well ordered $T$ in the following way
[the full construction is given in theorem A.2 in the appendix]: 
partition $T$ into a disjoint union of sub-branches, 
indexed by the nodes of a well founded tree $\Ga$ and reduce the problem of a 
well ordering of \/ $T$ to a problem of a well ordering of \/ $\Ga$. 
At the first step we pick a branch of \/ $T$, call it $A_{\lan\,\ran}$ and 
represent $T$ as \/ $A_{\lan\,\ran}\cup\bpl_{\et\in\lan\,\ran^+}T_\et$
(where for $\ta\in\Ga$, $\ta^+$ is the set $\{\nu: \nu$ an immediate successor  
of $\ta$ in $\Ga\}$ ). At the second step we pick a branch $A_\et$ in each 
$T_\et$ and represent $T_\et$ as $A_\et\cup\bpl_{\nu\in\et^+}T_\nu$. 
By tameness we finish after $\om$ steps getting $T=\cup_{\et\in\Ga}A_\et$ 
and the well ordering of $T$ is induced by the lexicographical well ordering 
of $\Ga$ and the well ordering of each $A_\et$ (which is scattered of 
Hdeg$\,\le k^*$). We can choose a sequence of parameters $\bb K_0$ (with 
length depending on $n^*$ and $k^*$ only) and a set of representatives 
$K=\{u_\et\in A_\et:\et\in\Ga\}$ and using $\bb K_0$ we can define a binary 
relation $<^*$ on $K$ where $u_\et<^*u_\nu$ will hold exactly when 
$\et\seg\nu$ in $\Ga$, thus we can define the structure of $\Ga$ in $T$. 
The sequence $\bb K_0$ will also enable us to define $T_\et$ and $A_\et$ from 
the representative $u_\et$ and define a well ordering of each $A_\et$. 

\\Consequently, the order between two nodes $x,y\in T$ will be determined by 
the well order of the $A_\et$'s (if they belong to the same $A_\et$) or the 
well ordering of $\Ga$ (if they belong to different $A_\et$'s). The well 
ordering of the sets $\et^+$ for $\et\in\Ga$ (hence the lexicographical well 
ordering of the well founded tree $\Ga$) will be again defined using 
$\bb K_0$.

\\What we'll do here in order to uniformize $\vf(X,Y,\bb Q)$ is the following:
given $Y\sb T$ we will use the decomposition $T=\cup_{\et\in\Ga}A_\et$ and 
the fact that each $A_\et$ is a scattered chain with Log$(A_\et)<l^*$, (hence 
satisfies the uniformization property), to define a unique $X_\et\sb A_\et$.
This will be done in such a way that when we glue the parts letting 
$X^*=\cup_{\et\in\Ga}X_\et$ we will still get $T\mo\vf(X,Y,\bb Q)$.

\\We will use the set of representatives $K$ and the fact that $A_\et$ and 
$T_\et$ are defined from $u_\et$ but we won't always mention $\bb K_0$.  
We will also rely on the fact that $\Ga$ is well founded (in fact, we only 
need to know that $\Ga$ does not have a branch of order type $\ge\om+1$).

\\So let $Y\sb T$ and we want to define some $X^*=X^*(Y,\bb Q)\sb T$. 
The proof will go as follows: for each $\et\in\Ga$ we will define partitions 
$\bbp^1(Y,\bb Q)_\et$ and $\bbp^2(Y,\bb Q)_\et$ of 
$K_{\et^+}:=\{u_\nu:\nu\in\et^+\}$ then, using the composition theorem 3.14 
and similarly to the proof of proposition 6.1, we will define a notion of 
coherence and let $\bb R^1(Y,\bb Q)_\et$ and $\bb R^2(Y,\bb Q)_\et$ be a pair 
that is coherent with $\bbp^1(Y,\bb Q)_\et$ and $\bbp^2(Y,\bb Q)_\et$. 
The union $\bb R^1(Y,\bb Q)=\cup_{\et\in\Ga}\bb R^1(Y,\bb Q)_\et$ is a 
partition of $K$ and $\th^n(A_\et;X_\et,Y\cap A_\et,\bb Q\cap A_\et)$ will be 
determined by the unique member of $\bb R^1(Y,\bb Q)$ to which $u_\et$ 
belongs. Moreover, we will be able to choose $X_\et$ uniquely and by coherence
$X^*=\cup_{\et\in\Ga}X_\et$ will satisfy $\vf(X,Y,\bb Q)$.

\line{\hfill{$\to$ [3.12.]}}

\\To get started let $T=A_{\lan\,\ran}\cup\bpl_{\et\in\lan\,\ran^+}T_\et$. 
Now as in definition 3.12 $K_{\lan\,\ran^+}$ has a natural structure of a 
chain with Log$(K_{\lan\,\ran^+})$=Log$(A_\lan\,\ran)<l^*$ and by theorem 
3.14(2) there is some $m=m(n,l)$ \st when $X\sb T$ is given, from 
$\th^m(A_{\lan\,\ran};X,Y,\bb Q)$ and 
$\lan\th^m(T_\et;X,Y,\bb Q):\et\in\lan\,\ran^+\ran$ we can compute 
$\th^n(T;X,Y,\bb Q)$. 

\\Let $\lan s_0,\ldots,s_{b-1}\ran$ be an enumeration of the 
the theories in $\T_{n+1,l+1}$. 

\\Define $\bbp^1(Y,\bb Q)_{\lan\,\ran}=
\lan P^1_0(Y,\bb Q)_{\lan\,\ran},\ldots,P^1_{b-1}(Y,\bb Q)_{\lan\,\ran}\ran$ 
a partition of $K_{\lan\,\ran^+}$ by 
$$\et\in\P^1_i(Y,\bb Q)_{\lan\,\ran} \iff \th^{n+1}(T_\et;Y,\bb Q)=s_i$$
By the previous remarks $\bbp^1(Y,\bb Q)_{\lan\,\ran}$ is definable from 
$u_{\lan\,\ran},K,Y,\bb Q$ (and $\bb K_0$).

\\Define $\bbp^2(Y,\bb Q)_{\lan\,\ran}=
\lan P^2_0(Y,\bb Q)_{\lan\,\ran},\ldots,P^2_{b-1}(Y,\bb Q)_{\lan\,\ran}\ran$ 
a partition of $K_{\lan\,\ran^+}$ by 
$$\et\in\P^1_i(Y,\bb Q)_{\lan\,\ran} \iff \th^{n+1}(A_\et;Y,\bb Q)=s_i$$
Again, $\bbp^2(Y,\bb Q)_{\lan\,\ran}$ is definable from 
$u_{\lan\,\ran},K,Y,\bb Q$ and $\bb K_0$.

\\Let $\lan t_0,\ldots,t_{a-1}\ran$ be an enumeration of the the theories in 
$\T_{n,l+2}$. 

\\A partition of $K_{\lan\,\ran^+}$, \/ 
$\bb R^1=\lan R^1_0,\ldots,R^1_{a-1}\ran$ is coherent with 
$\bbp^1(Y,\bb Q)_{\lan\,\ran}$ if 
$P^1_i(Y,\bb Q)_{\lan\,\ran}\cap R^1_j\ne\em$ implies 

\\``for every tree $S$ and $B,\bb C\sb S$ with $\lg(\bb C)=l$, 
if $\th^{n+1}(S;B,\bb C)=s_i$ then there is $A\sb S$ \st 
$\th^n(S;A,B,\bb C)=t_j$''.

\\Similarly a partition of $K_{\lan\,\ran^+}$, \/ 
$\bb R^2=\lan R^2_0,\ldots,R^2_{a-1}\ran$ is coherent with 
$\bbp^2(Y,\bb Q)_{\lan\,\ran}$ if 
$P^2_i(Y,\bb Q)_{\lan\,\ran}\cap R^2_j\ne\em$ implies 

\\``for every chain $S$ and $B,\bb C\sb S$ with $\lg(\bb C)=l$, 
if $\th^{n+1}(S;B,\bb C)=s_i$ then there is $A\sb S$ \st 
$\th^n(S;A,B,\bb C)=t_j$''.

\\Finaly, a pair of partitions of $K_{\lan\,\ran^+}$, \/ 
$\lan\bb R^1,\bb R^2\ran$ is \underbar{$t^*$-coherent} with the pair 
$\lan\bbp^1(Y,\bb Q)_{\lan\,\ran},\bbp^2(Y,\bb Q)_{\lan\,\ran}\ran$ if 

\\(1) $\bb R^1$ is coherent with $\bbp^1(Y,\bb Q)_{\lan\,\ran}$, 

\\(2) $\bb R^2$ is coherent with $\bbp^2(Y,\bb Q)_{\lan\,\ran}$, \ and

\\(3) For every $X\sb T$, \ if $\th^n(A_{\lan\,\ran};X,Y,\bb Q)=t^*$ and if
for every $\et\in\lan\,\ran^+$ \/ 
$\big[\th^n(T_\et;X,Y,\bb Q)=t_i\iff u_\et\in R^1_i\big]$, then 
$T\mo\vf(X,Y,\bb Q)$.

\\As $T\mo(\ex X)\vf(X,Y,\bb Q)$ there are $t^*$ (that will be fixed from now 
on), $\bb R^1$ and $\bb R^2$ \st 
$\lan\bb R^1,\bb R^2\ran$ is $t^*$-coherent with the pair 
$\lan\bbp^1(Y,\bb Q)_{\lan\,\ran},\bbp^2(Y,\bb Q)_{\lan\,\ran}\ran$. 

\\Moreover, ``$\lan\bb R^1,\bb R^2\ran$ is $t^*$-coherent with the 
pair $\lan\bbp^1(Y,\bb Q)_{\lan\,\ran},\bbp^2(Y,\bb Q)_{\lan\,\ran}\ran$'' 

\\is determined by $\th^k(K_{\lan\,\ran^+};\bb R^1,\bb R^2,
\bbp^1(Y,\bb Q)_{\lan\,\ran},\bbp^2(Y,\bb Q)_{\lan\,\ran})$ where $k$ depends 
only on $n$ and $l$. 

\\The first two clauses are clear (since $a$ and $b$ are finite) and for the 
third clause use theorem 3.14(2). So the statement is expressed by a p.u 
formula $\psi^1(\bb R^1,\bb R^2,
\bbp^1(Y,\bb Q)_{\lan\,\ran},\bbp^2(Y,\bb Q)_{\lan\,\ran})$ of depth $k$.

\\As by a previous remark Log$(K_{\lan\,\ran^+})<l^*$ there is 
$\bb S_{\lan\,\ran}\sb K_{\lan\,\ran^+}$ and a formula 
$\psi_{\lan\,\ran}(\bb U_1,\bb U_2,\bb W_1,\bb W_2,\bb S_{\lan\,\ran})$ 
that uniformizes $\psi^1$. 

\\To conclude the first step use Log$(A_{\lan\,\ran})<l^*$ to 
define, by a formula 
$\theta_{\lan\,\ran}(X,Y\cap A_{\lan\,\ran},\bb Q\cap A_{\lan\,\ran},\bb O_{\lan\,\ran})$ 
and a sequence of parameters $\bb O_{\lan\,\ran}\sb A_{\lan\,\ran}$, 
a unique $X_{\lan\,\ran}\sb A_{\lan\,\ran}$ that will satisfy 
$\th^n(A_\et;X_{\lan\,\ran},Y,\bb Q)=t^*$. 

\\

\\The result of the first step is the following: 

\\a) we have defined $X_{\lan\,\ran}\sb A_{\lan\,\ran}$ using 
$\bb O_{\lan\,\ran}\sb A_{\lan\,\ran}$ and $\theta_{\lan\,\ran}$. 
$X_{\lan\,\ran}$ is the intesection of the eventual $X^*$ with 
$A_{\lan\,\ran}$.

\\b) we have chosen $\bb R^1_{\lan\,\ran^+},\bb R^2_{\lan\,\ran^+}\sb 
K_{\lan\,\ran^+}$ using $\psi$ and $\bb S_{\lan\,\ran}$. 

\\c) $\bb R^1_{\lan\,\ran^+}$ and $\bb R^2_{\lan\,\ran^+}$ tell us what are
(for $\et\in\lan\,\ran^+$) the theories
$\th^n(T_\et;X^*,Y,\bb Q)$ and $\th^n(A_\et;X_\et,Y,\bb Q)$ respectively: 
if $u_\et\in R^1_i$ then the eventual $X^*\cap T_\et\sb T_\et$ 
will satisfy $\th^n(T_\et;X^*\cap T_\et,Y,\bb Q)=t_i$ and if 
$u_\et\in R^2_j$ then then the soon to be defined $X_\et\sb A_\et$ will 
satisfy $\th^n(A_\et;X_\et,Y,\bb Q)=t_j$.

\\We will proceed by induction on the level of $\et$ in $\Ga$ (remember, all 
the levels are $<\om$) to define $\bb S_\et,\bb O_\et\sb A_\et$ and 
$\bb R^1_{\et^+},\bb R^2_{\et^+}\sb K_{\et^+}$ and $X_\et\sb T_\et$.

\\The induction step: 

\\We are at $\nu\in\Ga$ where $\nu\in\et^+$ and we want to define 
$\bb S_\nu,\bb O_\nu\sb A_\nu$, \/ 
$\bb R^1_{\nu^+},\bb R^2_{\nu^+}\sb K_{\nu^+}$ and $X_\nu\sb T_\nu$.
Now as $\bb R^1_{\et^+}$ and $\bb R^2_{\et^+}$ are defined, 
$u_\nu$ belongs to one member of $\bb R^1_{\et^+}$ say the $i_1$'th 
and to one member of $\bb R^2_{\et^+}$ say the $i_2$'th. This implies that 
there is some $X'_\nu\sb T_\nu$ \st $\th^n(T_\nu;X'_\nu,Y,\bb Q)=t_{i_1}$ and 
$\th^n(A_\nu;X'_\nu\cap A_\nu,Y,\bb Q)=t_{i_2}$. 

\\Let $\bbp^1(Y,\bb Q)_\nu$ and $\bbp^2(Y,\bb Q)_\nu$ be partitions of 
$K_{\nu^+}$ that are defined as in the first step by saying, 
for $\ta\in\nu^+$, what are $\th^{n+1}(T_\ta;Y,\bb Q)$ and 
$\th^{n+1}(A_\ta;Y,\bb Q)$. 
$\lan\bb R^1_{\nu^+},\bb R^2_{\nu^+}\ran\sb K_{\nu^+}$ will be a pair that is 
\underbar{$t_{i_1},t_{i_2}$-coherent} with 
$\lan\bbp^1(Y,\bb Q)_\nu,\bbp^2(Y,\bb Q)_\nu\ran$ that is:

\\(1) $\bb R^1_{\nu^+}$ is coherent with $\bbp^1(Y,\bb Q)_\nu$, 

\\(2) $\bb R^2_{\nu^+}$ is coherent with $\bbp^2(Y,\bb Q)_\nu$, \ and

\\(3) For every $X\sb T_\nu$ if $\th^n(A_\nu;X,Y,\bb Q)=t_{i_2}$ and 
for every $\ta\in\nu^+$ \/ 
$\big[\th^n(T_\ta;X,Y,\bb Q)=t_i\iff u_\ta\in$ the $i$'th member of 
$R^1_{\nu^+}\big]$, then $\th^n(T_\nu;X,Y,\bb Q)=t_{i_1}$.

\\Using Log$(K_{\nu^+})<l^*$ choose $\bb S_{\nu^+}\sb K_{\nu^+}$ and
$\psi_{i_1,i_2}(\bb R^1,\bb R^2,\bbp^1(Y,\bb Q)_{\nu^+},
\bbp^2(Y,\bb Q)_{\nu^+},\bb S_{\nu^+})$ that uniformizes the formula that says 
``$\lan\bb R^1,\bb R^2\ran$ is $t_{i_1},t_{i_2}$-coherent with 
$\lan\bbp^1(Y,\bb Q)_\nu,\bbp^2(Y,\bb Q)_\nu\ran$''. We may assume that 
$\psi_{i_1,i_2}$ depends only on $i_1$ and $i_2$ and that $\lg(\bb S_{\nu^+})$ 
is constant. 

\\Use Log$(A_\nu)<l^*$ to define, by a formula 
$\theta_{i_2}(X,Y\cap A_\nu,\bb Q\cap A_\nu,\bb O_\nu)$ and a sequence of 
parameters $\bb O_\nu\sb A_\nu$, a unique $X_\nu\sb A_\nu$ that will satisfy 
$\th^n(A_\nu;X_\nu,Y,\bb Q)=t_{i_2}$. Again, we may assume that 
$\theta_{i_2}$ depends only on $i_2$ and that $\lg(\bb O_\nu)$ is constant. 

\\So $\bb S_\nu$, $\bb O_\nu$, $\bb R^1_{\nu^+}$, $\bb R^2_{\nu^+}$ 
and $X_\nu$ are defined and we have concluded the inductive step.
(Note that nothing will realy go wrong if $\nu$ doesn't have any successors 
in $\Ga$).

\\

\\Let $\bb O=\cup_{\et\in\Ga}\bb O_\et$, \/ $\bb S=\cup_{\et\in\Ga}\bb S_\et$. 
The uniformizing formula $U(X,Y,\bb Q,\bb O,\bb S,K,\bb K_0)$ says:

\\``$X\cap A_{\lan\,\ran}$ is defined as in the first step, \ and

\\for every pair of partitions $\lan\bbp^1,\bbp^2\ran$ of $K$ that agrees
on each $K_{\et^+}$ with [the definable] 

\\$\lan\bbp^1_{\et^+}(Y,\bb Q),\bbp^2_{\et^+}(Y,\bb Q)\ran$,  
(and agrees with $\lan\bbp^1_{\lan\,\ran},\bbp^2_{\lan\,\ran}\ran$ on 
$K_{\lan\,\ran^+}$), and

\\for every $\lan\bb R^1,\bb R^2\ran$ that is a [in fact the only] pair of 
partitions that satisfies for every $u_\et\in K$: 
if $u_\et\in P^1_{i_1}\cap P^2_{i_2}$ then
$\psi_{i_1,i_2}(\bb R^1\cap K_{\et^+},\bb R^2\cap K_{\et^+},
\bbp^1\cap K_{\et^+},\bbp^2\cap K_{\et^+},\bb S\cap K_{\et^+})$ holds,
(and agrees with $\lan\bb R^1_{\lan\,\ran},\bb R^2_{\lan\,\ran}\ran$ on 
$K_{\lan\,\ran^+}$),

\\for every $u_\et\in K$ if $u_\et\in R^2_i$ then 
$\theta_i(X\cap A_\et,Y\cap A_\et,\bb Q\cap A_\et,\bb O\cap A_\et)$ holds.''

\\

\\$U(X,Y,\bb Q,\bb O,\bb S,K,\bb K_0)$ does the job because it 
defines $X\cap A_\et$ uniquely on each $A_\et$ and because, (by the conditions 
of coherence) the union of the parts, $X$, satisfies $\vf(X,Y,\bb Q)$. Note 
also that $U$ does not depend on $Y$.
\vv \qed
\sec
%
%
%
%
%
%
\hh 7. Hopelessness of General Partial Orders   \par
%
%
%
\pt Theorem 7.1. Every partial order $P$ can be embedded in a partial order 
$Q$ in which $P$ is first-order-definably well orderable.  
\pd Proof.
\vv \qed
\sec
%
%
%
%
%
%
\hh Appendix  \par
%
%
%
\pt Lemma A.1. Let $C$ be a scattered chain with ${\rm Hdeg}(C)=n$. Then there are 
$\bb P\sb C$, $\lg(\bb P) = n-1$, and a formula (depending on $n$ only)  
$\varphi_n(x,y,\bb P)$ that defines a well ordering of\/ $C$.  \vv
\pd Proof. By induction on $n={\rm Hdeg}(C)$: 

\\$n\le 1$: \ ${\rm Hdeg}(C)\le 1$ implies $(C,<_C)$ is well ordered or inversely 
well ordered. A well ordering of\/ $C$ is easily definable from $<_C$.

\\${\rm Hdeg}(C)=n+1$: Suppose $C=\sum_{i\in I}C_i$ and each 
$C_i$ is of Hausdorff degree $n$. By the induction hypothesis there are a 
formula $\varphi_n(x,y,\bb Z)$ and a sequence $\lan \bb P^i: i\in I \ran$\  
with $\bb P^i\sb C_i$, $\bb P^i = \lan P_1^i,\ldots,P_{n-1}^i \ran$ \st
$\varphi_n(x,y,\bb P^i)$ defines a well ordering of\/ $C_i$. 

\\Let for $0<k<n$, $P_k := \cup_{i\in I}P^i_k$  (we may assume that the union 
is disjoint) and $P_n := \cup\{C_i : i {\rm \ even}\}$.

\\We will define an equivalence relation $\sim$ by $x\sim y$ iff 
$\bigwedge_i (x\in C_i \Leftrightarrow y\in C_i)$. 

\\$\sim$ and $[x]$, (the equivalence class of an element $x$), are easily 
definable from $P_n$ and $<_C$. We can also decide from $P_n$ if\/ $I$ is well 
or inversely well ordered (by looking at subsets of\/ $C$ consisted of 
nonequivalent elements) and define $<'$ to be $<$ if\/ $I$ is well ordered 
and the inverse of\/ $<$ if not.
$\varphi_{n+1}(x,y,P_1,\ldots,P_n)$ will be 
defined by: 
$$\varphi_{n+1}(x,y,\bb P) \Leftrightarrow 
\big[ x\not\sim y \ \&\  x<'y\big] \vee 
\big[x\sim y \ \&\ \varphi_{n}(x,y,P_1\cap[x],\ldots,P_{n-1}\cap[x])\big]$$ 

\\$\varphi_{n+1}(x,y,\bb P)$ well orders $C$. 
\vv \qed
\pt Theorem A.2. Let $T$ be a tame tree. If\/ $\bina$ is not embeddable in 
$T$ then there are $\bb Q\sb T$ and a monadic formula $\varphi(x,y,\bb Q)$ 
that defines a well ordering of\/ $T$. \vv
\pd Proof. Assume $T$ is $(n^*,k^*)$ tame, recall definitions 4.1 and 4.2 and 
remember that for every $x \in T$, $rk(x)$ is well defined (i.e. $<\infty$).
We will partition $T$ into a disjoint union of sub-branches, 
indexed by the nodes of a well founded tree $\Ga$ and reduce the problem of a 
well ordering of\/ $T$ to a problem of a well ordering of\/ $\Ga$.

Step 1. Define by induction on $\al$ a set $\Ga_\al\sb\thinspace^\al Ord$ 
(this is a our set of indices), for every $\et\in\Ga_\al$ define a tree 
$T_\et\sb T$ and a branch $A_\et\sb T_\et$.

\\$\al=0$ : $\Ga_0$ is $\{\lan \ran\}$, $T_{\lan \ran}$ is $T$ and 
$A_{\lan \ran}$ is a branch (i.e. a maximal linearily ordered subset) of\/ $T$.

\\$\al=1$ : Look at $(T\setminus A_{\lan \ran})/\sim^1_{A_{\lan \ran}}$, it's 
a disjoint union of trees and name it $\lan T_{\lan i\ran}:i<i^* \ran$,
let $\Ga_1:=\{\lan i\ran:i<i^*\}$ and for every $\lan i\ran\in\Ga_1$ let 
$A_{\lan i\ran}$ be a branch of\/ $T_{\lan i\ran}$.

\\$\al=\be+1$ : For $\et\in\Ga_\be$ denote 
$(T_\et\setminus A_\et)/\sim^1_{A_\et}$ by $\{T_{\ee\et,i}:i<i_\et\}$, let  
$\Ga_\al=\{\ee\et,i:\et\in\Ga_\be,\ i<i_\et\}$ and choose $A_{\ee\et,i}$ to 
be a branch of\/ $A_{\ee\et,i}$.

\\$\al$ limit: Let $\Ga_\al=
\{\et\in\thinspace^\al Ord : \wedge_{\be<\al}\et\res_\be\in\Ga_\be,\ 
\wedge_{\be<\al}T_{\et\res_\be}\not=\emptyset\}$, 
let for $\et\in\Ga_\al$ $T_\et=\cap_{\be<\al}T_{\et\res_\be}$ and 
$A_\et$ a branch of\/ $T_\et$. ($T_\et$ may be empty).

\\Now, at some stage $\al\le|T|^+$ we have $\Ga_\al=\emptyset$ and let 
$\Ga=\cup_{\be<\al}\Ga_\be$. Clearly $\{A_\et:\et\in\Ga\}$ is a partition of 
$T$ into disjoint sub-branches.

\\Notation: having two trees $T$ and $\Ga$, to avoid confusion, we use 
$x,y,s,t$ for nodes of\/ $T$ and $\et,\nu,\si$ for nodes of\/ $\Ga$.

Step 2. We want to show that $\Ga_\om=\emptyset$ hence $\Ga$ is a well 
founded tree. Note that we made no restrictions on the choice of 
the $A_\et$'s and we add one now in order to make the above statement true.
Let $\ee\et,i\in\Ga$ define $A_{\et,i}$ to be the sub-branch 
$\{t\in A_\et:(\all s\in A_{\ee\et,i})[rk(t)\le rk(s)]\}$ and 
$\ga_{\et,i}$ to be $rk(t)$ for some $t\in A_{\et,i}$. By 5.5(1) and the 
inexistence of a stricly decreasing sequence of ordinals, 
$A_{\et,i}\not=\emptyset$ and $\ga_{\et,i}$ is well defined. Note also that 
$s\in A_{\ee\et,i} \imp rk(s)\le\ga_{\et,i}$. 

\\\underbar {Proviso}: For every $\et\in\Ga$ and $i<i_\et$ the sub-branch 
$A_{\ee\et,i}$ contains every $s\in T_{\ee\et,i}$ with $rk(s)=\ga_{\et,i}$. 

\\Following this we claim: ``$\Ga$ does not contain an infinite, stricly 
increasing sequence''.  Otherwise let $\{\et_i\}_{i<\om}$ be one, and 
choose $s_n\in A_{\et_n,\et_{n+1}(n)}$ (so $s_n\in A_{\et_n}$). Clearly 
$rk(s_n)\ge rk(s_{n+1})$ and by the proviso we get 
$$rk(s_n)=rk(s_{n+1})\imp rk(s_{n+1})>rk(s_{n+2})$$
\\therefore $\{rk(s_n)\}_{n<\om}$ contains an infinite, stricly decreasing 
sequence of ordinals which is absurd.

Step 3. Next we want to make ``$x$ and $y$ belong to the same $A_\et$'' 
definable. 

\\For each $\et\in\Ga$ choose $s_\et\in A_\et$, and let $Q\sb T$ be the set
of representatives.
Let $h\colon T\to\{d_0,\ldots,d_{n^*-1}\}$ be a colouring that satisfies:
$h\res_{A_{\lan \ran}}=d_0$ and for every $\ee\et,i\in\Ga$, 
$h\res_{A_{\ee\et,i}}$ is constant and, when $j<i$ and 
$s_{\ee\et,j}\sim^0_{A_\et}s_{\ee\et,j}$ we have  
$h\res_{A_{\ee\et,i}}\not=h\res_{A_{\ee\et,j}}$. This can be done as $T$ is 
$(n^*,d^*)$ tame.

\\Using the parameters $D_0,\ldots,D_{n^*-1}$ ($x\in D_i$ iff\/ $h(x)=d_i$),
we can define $\vee_\et x,y\in A_\et$ by ``$x,y$ are comparable and the 
sub-branch $[x,y]$ (or $[y,x]$) has a constant colour''.

Step 4. As every $A_\et$ has Hausdorff degree at most $k^*$, we can define a 
well ordering of it using parameters $P^\et_1,\ldots,P^\et_{k^*}$ and by 
taking $\bb P$ to be the (disjoint) union of the $\bb P^\et$'s we can define 
a partial ordering on $T$ which well orders every $A_\et$.

\\By our construction $\et\seg\nu$ if and only if there is an element in 
$A_\nu$ that `breaks' $A_\et$ i.e. is above a proper initial segment of 
$A_\et$. (Caution, if\/ $T$ does not have a root this may not be the case for  
$\lan \ran$ and a $<n^*$ number of\/ $\lan i\ran$'s and we may need parameters 
for expressing that). Therefore, as by step 3 ``being in the same $A_\et$'' 
is definable, we can define a partial order on the sub-branches $A_\et$ 
(or the representatives $s_\et$) by $\et\seg\nu\imp A_\et\le A_\nu$.  

\\Next, note that ``$\nu$ is an immediate successor of\/ $\et$ in $\Ga$'' is 
definable as a relation between $s_\nu$ and $s_\et$ hence the set 
$A_\et^+ := A_\et\cup\{s_{\ee\et,i}\}$ is definable from $s_\et$.  Now 
the order on $A_\et$ induces an order on $\{s_{\ee\et,i}/\sim^0_{A_\et}\}$ 
which is can be embedded in the complition of\/ $A_\et$ hence has Hdeg$\le k^*$. 
Using additional parameters $Q^\et_1,\ldots,Q^\et_{k^*}$, we have a definable 
well ordering on $\{s_{\ee\et,i}/\sim^0_{A_\et}\}$. As for the ordering on 
each $\sim^1_{A_\et}$ equivalence class (finite with $\le n^*$ elements),
define it by their colours (i.e. the element with the smaller colour is the 
smaller according to the order).

\\Using $\bb D$, $\bb P$, $Q$ and $\bb Q=\cup_\et\bb Q^\et$ we can define a 
partial ordering which well orders each $A_\et^+$ in such a way that every
$x\in A_\et$ is smaller then every $s_{\ee\et,i}$.

\\Summing up we can define (using the above parameters) a partial order 
on subsets of\/ $T$ that well orders each $A_\et$, orders sub-branches 
$A_\et$, $A_\nu$ when the indices are comparable in $\Ga$ and well 
orders all the ``immediate successors'' sub-branches of a sub-branch $A_\et$.

Step 5.  The well ordering of\/ $T$ will be defined by $x<y \iff$

\\a) $x$ and $y$ belong to the same $A_\et$ and $x<y$ by the well order on 
$A_\et$;  \ or

\\b) $x\in A_\et$, $y\in A_\nu$ and $\et\seg\nu$;  \ or

\\c) $x\in A_\et$, $y\in A_\nu$,  $\si=\et\wedge\nu$ in $\Ga$ (defined as a 
relation between sub-branches), $\ee\si,i\seg\et$, $\ee\si,j\seg\nu$  
and $s_{\ee\si,i}<s_{\ee\si,j}$ in the order of\/ $A^+_\si$.

\\Note, that $<$ is a linear order on $T$ and every $A_\et$ is a convex 
and well ordered sub-chain. Moreover $<$ is a linear order on $\Ga$ and the 
order on the $s_\et$'s is isomorphic to a lexicographic order on $\Ga$.

\\Why is the above (which is clearly definable with our parameters) a well 
order?  Because of the above note and because a lexicographic ordering of a  
well founded tree is a well order, provided that immediate successors are 
well ordered. In detail, assume $X=\{x_i\}_{i<\om}$ is a stricly decreasing 
sequence of elements of\/ $T$.  Let $\et_i$ be the unique node in $\Ga$ \st 
$x_i\in A_{\et_i}$ and by the above note w.l.o.g 
$i\not=j\imp\et_i\not=\et_j$.  By the well foundedness of\/ $\Ga$ and clause 
(b) we may also assume w.l.o.g that the $\et_i$'s form an anti-chain in 
$\Ga$. Look at  $\nu_i:=\et_1\wedge\et_i$  which is constant for infinitely 
many $i$'s and w.l.o.g equals to $\nu$ for every $i$.  Ask:

\\$(*)$ is there is an infinite $B\sb\om$ \st 
$i,j\in B\imp x_i\sim^0_{A_\nu}x_j$\ ?

\\If this occurs we have $\nu_1\not=\nu$ with $\nu\seg\nu_1$ \st for 
some infinite $B'\sb B\sb\om$ we have $i\in B'\imp \nu_1\seg\et_i$. 
(use the fact that $\sim^1_{A_\nu}$ is finite).
W.l.o.g $B'=\om$ and we may ask if $(*)$ holds for $\nu_1$. 
Eventually, since $\Ga$ does not have an infinite branch, we will have a 
negative answer to $(*)$. We can conclude that w.l.o.g there is 
$\nu\in\Ga$ \st $i\not=j\imp x_i\not\sim^0_{A_\nu}x_j$ i.e. the $x_i$'s 
``break'' $A_\nu$ in ``different places''.  

\\Define now $\nu_i$ to be the unique immediate successor of\/ $\nu$ \st 
$\nu_i\seg\et_i$. The set $S=\{s_{\nu_i}\}_{i<\om }\sb A^+_\nu$ is well 
ordered by the well ordering on $A^+_\nu$ and by clause (c) in the definition 
of\/ $<$,  $x_i>x_j\iff\nu_i>\nu_j$ so $S$ is an infinite stricly decreasing 
subset of\/ $A^+_\nu$ -- a contradiction.

\\This finishes the proof that there is a definable well order of\/ $T$.
\vv \qed

\end